\documentclass[11pt]{article}
\usepackage[margin=0.85in]{geometry}
\usepackage{authblk}
\usepackage[style=numeric, backend=biber, bibencoding=utf8]{biblatex}
\usepackage{hyperref}
\usepackage{graphicx}
\usepackage{subcaption}
\usepackage{wrapfig}
\usepackage{soul}

\usepackage{amsmath, amssymb, amsbsy}

\newcommand{\reals}{\mathbb{R}}

\newcommand{\realspace}[1]{\mathbb{R}^{#1}}

\newcommand{\deriv}[2]{\frac{d #1}{d #2}}
\newcommand{\pderiv}[2]{\frac{\partial #1}{\partial #2}}
\newcommand{\Hone}{H^1}
\newcommand{\Identity}{\boldsymbol{I}}
\newcommand{\hadamard}[2]{{#1}\odot{#2}}

\newcommand{\fulldomain}{\Omega}
\newcommand{\femspace}{X}
\newcommand{\bilin}{a}
\newcommand{\lin}{f}
\newcommand{\femsol}{u}
\newcommand{\testfun}{v}

\newcommand{\paramspace}{M}
\newcommand{\numparams}{n^m}
\newcommand{\prscparam}{\boldsymbol{\mu}}

\newcommand{\domaink}[1]{\Omega_{#1}}
\newcommand{\domaini}{\domaink{i}}

\newcommand{\numparamsi}{n^m_i}
\newcommand{\prscparamk}[1]{\prscparam_{#1}}
\newcommand{\prscparami}{\prscparamk{i}}
\newcommand{\bilink}[1]{\bilin_{#1}}
\newcommand{\bilini}{\bilink{i}}
\newcommand{\femspacek}[1]{\femspace_{#1}}
\newcommand{\femspacei}{\femspacek{i}}
\newcommand{\ncomps}{n_c}
\newcommand{\testfunk}[1]{\testfun_{#1}}
\newcommand{\femsolk}[1]{\femsol_{#1}}
\newcommand{\testfuni}{\testfunk{i}}
\newcommand{\femsoli}{\femsolk{i}}
\newcommand{\rfemsol}{\tilde{\femsol}}

\newcommand{\portkk}[2]{\gamma_{{#1},{#2}}}
\newcommand{\portik}[1]{\portkk{i}{#1}}
\newcommand{\portij}{\portkk{i}{j}}
\newcommand{\nportsk}[1]{n_{#1}^\gamma}
\newcommand{\nportsi}{\nportsk{i}}
\newcommand{\nportdofsp}{n^\Gamma_p}
\newcommand{\indexset}{\pi_p}
\newcommand{\globportp}{\Gamma_p}
\newcommand{\nglobports}{n^\Gamma}
\newcommand{\portmapi}{\mathcal{G}_i}
\newcommand{\portspaceij}{P_{i,j}}
\newcommand{\portspacekk}[2]{P_{{#1},{#2}}}
\newcommand{\nportdofsij}{n^\gamma_{i,j}}
\newcommand{\pbasisijk}{\chi_{i,j,k}}
\newcommand{\lpbasisijk}{\psi_{i,j,k}}
\newcommand{\lpbasiskkk}[2]{\psi_{{#1}, {#2}, k}}
\newcommand{\bubspacei}{B_i}
\newcommand{\gpbasiskk}[1]{\Psi_{{#1}, k}}
\newcommand{\gpbasispk}{\gpbasiskk{p}}
\newcommand{\bubijk}{b_{i,j,k}}
\newcommand{\bubfi}{b^f_i}
\newcommand{\prscsolkk}[2]{\mathbb{U}_{{#1},{#2}}}

\newcommand{\ifunckkk}[3]{\phi_{{#1},{#2},{#3}}}
\newcommand{\ifuncijk}{\ifunckkk{i}{j}{k}}
\newcommand{\Ifunckk}[2]{\Phi_{{#1}, {#2}}}
\newcommand{\Ifuncpk}{\Ifunckk{p}{k}}
\newcommand{\scstiff}{\mathbb{K}}
\newcommand{\scforce}{\mathbb{F}}
\newcommand{\scsol}{\mathbb{U}}
\newcommand{\scstiffk}[1]{\scstiff_{#1}}
\newcommand{\sccompstiffik}[1]{\mathbb{K}^i_{#1}}
\newcommand{\scforcek}[1]{\scforce_{#1}}
\newcommand{\rpbasisijk}{\tilde{\chi}_{i,j,k}}
\newcommand{\nrportdofsij}{\tilde{n}^\gamma_{i,j}}

\newcommand{\rscstiff}{\tilde{\scstiff}}
\newcommand{\rscforce}{\tilde{\scforce}}
\newcommand{\rscsol}{\tilde{\scsol}}

\newcommand{\densityvec}{\boldsymbol{\rho}}
\newcommand{\densityk}[1]{\rho_{#1}}
\newcommand{\densityi}{\densityk{i}}
\newcommand{\densitymin}{\densityk{min}}

\newcommand{\simpscaling}{s}
\newcommand{\elastictensor}{\mathbb{C}}
\newcommand{\lamemu}{\mu}
\newcommand{\lamelambda}{\lambda}
\newcommand{\force}{F}
\newcommand{\youngsi}{E_i}

\newcommand{\objective}{\mathcal{J}}
\newcommand{\stress}{\sigma}
\newcommand{\stressmax}{\sigma_{max}}
\newcommand{\heurstressmax}{\hat{\sigma}_{max}}
\newcommand{\vonmises}{\sigma_{vm}}
\newcommand{\relaxedstress}{\sigma_r}
\newcommand{\naggs}{n_{agg}}

\newcommand{\constrainti}{g_m}
\newcommand{\aggdomaini}{\Omega^{agg}_m}
\newcommand{\stressvecixx}{\boldsymbol{\sigma}^m_{xx}}
\newcommand{\stressveciyy}{\boldsymbol{\sigma}^m_{yy}}
\newcommand{\stressvecixy}{\boldsymbol{\sigma}^m_{xy}}
\newcommand{\stressopixx}{S^m_{xx}}
\newcommand{\stressopiyy}{S^m_{yy}}
\newcommand{\stressopixy}{S^m_{xy}}
\newcommand{\ifuncxpx}{\tilde{\Phi}^m_{x,x}}
\newcommand{\ifuncypy}{\tilde{\Phi}^m_{y,y}}
\newcommand{\ifuncxpy}{\tilde{\Phi}^m_{x,y}}
\newcommand{\ifuncypx}{\tilde{\Phi}^m_{y,x}}
\newcommand{\vmquadi}{\boldsymbol{\sigma}^m_{VM}}
\newcommand{\rstressveci}{\boldsymbol{\sigma}_r^m}
\newcommand{\aggrhoi}{\boldsymbol{\rho}^m_{agg}}

\newcommand{\prolongi}{P_m}
\newcommand{\adjointsol}{\boldsymbol{\lambda}}

\usepackage{color}
\definecolor{Bronze}{rgb}{0.8,0.5,0.2}
\definecolor{darkred}{rgb}{0.55, 0.0, 0.0}

\newcommand{\keywords}[1]{\small
    \textbf{\textit{Keywords---}} #1
}

\title{Stress-constrained topology optimization of lattice-like structures
using component-wise reduced order models}

\author[1]{Sean McBane\footnote{Corresponding author; sean@oden.utexas.edu}}
\author[2]{Youngsoo Choi}
\author[1]{Karen Willcox}
\affil[1]{Oden Institute for Computational Engineering \& Sciences, University of Texas, Austin, TX 78712}
\affil[2]{Center for Applied Scientific Computing, Lawrence Livermore National Laboratory, Livermore, CA 94550}

\date{}

\begin{document}
\maketitle

\begin{abstract}
    Lattice-like structures can provide a combination of high stiffness
with light weight that is useful in many applications, but a resolved
finite element mesh of such structures results in a computationally
expensive discretization. This computational expense may be particularly
burdensome in many-query applications, such as optimization.
We develop a stress-constrained topology optimization method
for lattice-like structures that uses component-wise reduced order models
as a cheap surrogate, providing accurate computation of stress fields while greatly
reducing run time relative to a full order model. We demonstrate the
ability of our method to produce large reductions in mass while
respecting a constraint on the maximum stress in a pair of test problems.
The ROM methodology provides a speedup of about 150x in forward solves
compared to full order static condensation and provides a relative error of less than
5\% in the relaxed stress.
\end{abstract}


\section{Introduction}
Lattice-like structures are advantageous in many applications due to their
combination of high stiffness and light weight. Applications include
aerospace structures \cite{opgenoord_willcox_2019}, medical uses such as
prostheses or implants \cite{xiao_et_al_2020}, and the traditional use of
trusses as supports in structural engineering.
Interest in lattice-like structures has further increased as
advances in additive manufacturing enable their production.
Their analysis using a standard finite element
method (FEM) requires a high-dimensional discretization to capture the complex
geometry, especially when accurate computation of localized quantities
-- for example, stress -- is needed. Optimization is then expensive, requiring
many evaluations of this high-dimensional model. In this work, we present a
stress-constrained topology optimization (TO) formulation for lattice-like structures
that leverages component-wise reduced order models (ROMs) as a surrogate, providing
accurate computation of the stress field while greatly reducing computational cost
relative to a FEM analysis.

TO of lattice-like structures is often based on asymptotic homogenization.
Such techniques are compelling when the final design should be
composed of material with a periodic structure. These techniques may solve a
macroscopic design problem in which the optimization variables control the
geometry of unit cells throughout the domain \cite{bendsoe_kikuchi_1988};
a microscopic design problem, designing a unit cell with particular properties
as first seen in \cite{sigmund_1994}; or a combination of the two as
in \cite{zhang_wang_kang_2019}. Homogenization has been applied to stress-constrained
problems, for example in \cite{cheng_bai_to_2019}.
Other methods also assume a periodic structure,
for example, the generalized and high-fidelity generalized method of cells \cite{aboudi_2004}.
ROMs have been applied to this family of techniques
\cite{ricks_et_al_2018,xia_breitkopf_2014,xia_breitkopf_2015};
however, together with homogenization, these methods share limiting assumptions
on length scale and periodicity. The length scale of unit cells should be much
less than that of macroscopic features, and the
approximation may not be accurate if material does not satisfy the periodicity
assumption. In particular, the behavior of stress near the boundary of the
material is not well understood \cite{cheng_bai_to_2019}.

Several other existing techniques in TO also apply to design of lattice-like structures.
Ground structure approaches for truss optimization \cite{zhang_et_al_ground_structure_2017,fairclough_et_al_2022}
rely on a beam model for truss members and make the optimization variables the
cross-sectional areas of truss members. This approximation is invalid when
truss members become too thick. \citeauthor{deng_to_2020} in \cite{deng_to_2020}
develop a method that projects the ground structure onto a background finite
element mesh, eliminating the thin beam assumption. In a similar vein, the family of methods
related to moving morphable components (MMC) \cite{guo_zhang_zhong_2014,norato_bell_tortorelli_2015}
represent a structure as a collection of discrete geometric components whose
shape and position are controlled by the optimization parameters, and model the
resulting structure by projecting component geometries on a background mesh.
\citeauthor{zhang_gain_norato_2017} in \cite{zhang_gain_norato_2017} and
\citeauthor{zhang_li_zhou_du_li_guo_2018} in \cite{zhang_li_zhou_du_li_guo_2018}
apply such projection-based methods to solve stress-constrained TO problems.
The shared shortcoming of these approaches for our purposes is the requirement
of a sufficiently fine background mesh to capture small-scale structural
features accurately. Substructuring methods accelerate solution of a high-resolution
model by rewriting equations in terms of a reduced set of degrees of freedom,
alleviating the computational cost associated with a fine discretization.
In \cite{our_first_paper}, we apply substructuring-based
model order reduction in the form of port-reduced static condensation (PRSC)
\cite{prsc} to compliance minimization. \citeauthor{wu_xia_wang_shi_2019}
solve the same problem using a related approach in \cite{wu_xia_wang_shi_2019};
\cite{wu_fan_xiao_yu_2020} and \cite{koh_kim_yoon_2020} are additional examples of
substructuring-based model order reduction to TO. Finally, there is an extensive body
of work applying conventional (i.e., density-based or level set) TO methods to
stress-constrained problems, and these approaches could be used to design lattice-like
structures given a fine enough discretization. See, for example,
\cite{le_stressbased_2010,holmberg_stressbased_2013} for examples in density-based
optimization to which we refer in our approach, and
\cite{kambampati_chung_kim_2021,picelli_et_al_2018} for examples of level set approaches
to stress-based optimization.

In this paper we leverage PRSC to efficiently solve stress-constrained
TO problems, which are not addressed in the previous literature on substructuring for TO.
Our approach is both a ground structure and a substructuring method. We constrain
the design space by choosing a fixed arrangement of subdomains (``components'')
of which all possible designs are a subset. We apply PRSC to construct a surrogate
model parameterized by a component-wise density parameter, which penalizes stiffness
using the SIMP scheme \cite{simp}. Then, we minimize the mass
of the structure subject to a constraint on the maximum stress implemented using
the \textit{qp}-relaxation of \citeauthor{bruggi_2008} \cite{bruggi_2008}
to address the ``singularity problem'' \cite{rozvany_1996},
and stress aggregation in the form of the Kreisselmeir-Steinhauser (KS) functional
to convert the infinite-dimensional stress constraint to a small number of
differentiable constraints that approximate the max function. Constraints are
defined by aggregating over non-overlapping aggregation domains, as in
\cite{le_stressbased_2010,holmberg_stressbased_2013}.
Finally, after optimization we postprocess by removing components with
densities less than a prescribed minimum value. Due to the component-wise formulation
this postprocessing results in a well-defined geometry without additional
postprocessing steps, unlike the result from an element-wise density based TO.

The remainder of this paper is organized as follows. In Section \ref{sec:prsc},
we provide a brief overview of PRSC, introducing the notation required for
the rest of our discussion. Section \ref{sec:cw_formulation} describes our
optimization formulation: the material model, constraint formulation,
objective function, and postprocessing methodology.
Following this description, in Section \ref{sec:numerical_results} we
present mass minimization results for an L-bracket and a cantilever beam
geometry, along with studies of ROM accuracy and performance. Finally,
we give our conclusions.

\section{An overview of port-reduced static condensation}
\label{sec:prsc}
Port-reduced static condensation (PRSC) is developed in
\cite{prsc} and further discussed in \cite{scrbe_nonlin, smetana_2015, smetana_patera_2016, iapichino_et_al_2016, our_first_paper}
and others. The discussion below is a high-level overview of the technique in which
we seek to introduce only the concepts and notation needed to describe the
use of PRSC for component-wise TO and perform the required sensitivity analysis.

PRSC applies to solution of an elliptic PDE on a domain
$\fulldomain \subset \realspace{n}$ (where $n = 2 \text{ or } 3$ in general, but we
focus on $n = 2$ here) given in weak form by
\begin{equation}\label{eq:global_weak_form}
    \bilin(\femsol, \testfun; \prscparam) = \lin(\testfun; \prscparam),
    \ \forall \testfun \in \femspace
\end{equation}
Here, $\bilin: \Hone(\fulldomain) \times \Hone(\fulldomain) \to \reals$ is a
bilinear form, taken to be coercive so that the problem \eqref{eq:global_weak_form}
is well-posed; $\lin: \Hone(\fulldomain) \to \reals$ is a linear form;
$\prscparam \in \paramspace$ is a vector in
$\paramspace \subset \realspace{\numparams}$, parametrizing both the linear
and bilinear forms; $\femspace \subset \Hone(\fulldomain)$ is a
finite dimensional function space arising from a finite element discretization
of \eqref{eq:global_weak_form} and incorporating essential boundary conditions; and
$\femsol, \testfun \in \femspace$ are the solution to \eqref{eq:global_weak_form}
and a test function, both residing in $\femspace$.

PRSC reduces the number of degrees of freedom in \eqref{eq:global_weak_form} by
applying static condensation, before further reducing the problem dimension
via projection-based model reduction. $\fulldomain$ is decomposed
into $\ncomps$ subdomains $\domaini, \ i \in \{1, \ldots, \ncomps\}$, termed ``components''.
The bilinear and linear forms may now be decomposed as
\begin{equation}\label{eq:cw_bilinear_form}
    a(u, v; \prscparam) = \sum_{i=1}^{\ncomps}
    \bilini\left(\femsoli, \testfuni; \prscparami\right)
\end{equation}
\begin{equation}\label{eq:cw_linear_form}
    f(v; \prscparam) = \sum_{i=1}^{\ncomps}
    \lin_i\left(\testfuni; \prscparami\right)
\end{equation}
where $\femsoli$ and $\testfuni$ are the restrictions of $\femsol$ and $\testfun$
to $\domaini$; $\bilini$ and $\lin_i$ are the restrictions of $a$ and $f$ to act
on functions in
$\femspace|_{\domaini}$; and $\prscparami \subset \realspace{\numparamsi}$,
with $\numparamsi$ the dimension of the parameter space for component $i$,
is a parameter vector containing only those
parameters in $\prscparam$ to which $\bilini$ and $\lin_i$ are sensitive.

Each component is assigned a set of $\nportsi$ ``ports'',
$\left\{\portij, \ j \in 1, \ldots, \nportsi\right\}$, which are subsets of the
boundary of $\domaini$ where $\domaini$ may (but is not required to) interface with
another component $\domaink{i'}$. That is, if $\domaini \cap
\domaink{i'}$ is non-empty, it corresponds to port $\portij$ and $\portkk{i'}{j'}$
respectively on $\domaini$ and $\domaink{i'}$. Ports may also be unconnected to
any neighboring component, or have a Dirichlet boundary applied. Ports on the
same component are assumed to be disjoint.

\begin{figure}
    \centering
    \includegraphics[width=0.8\textwidth]{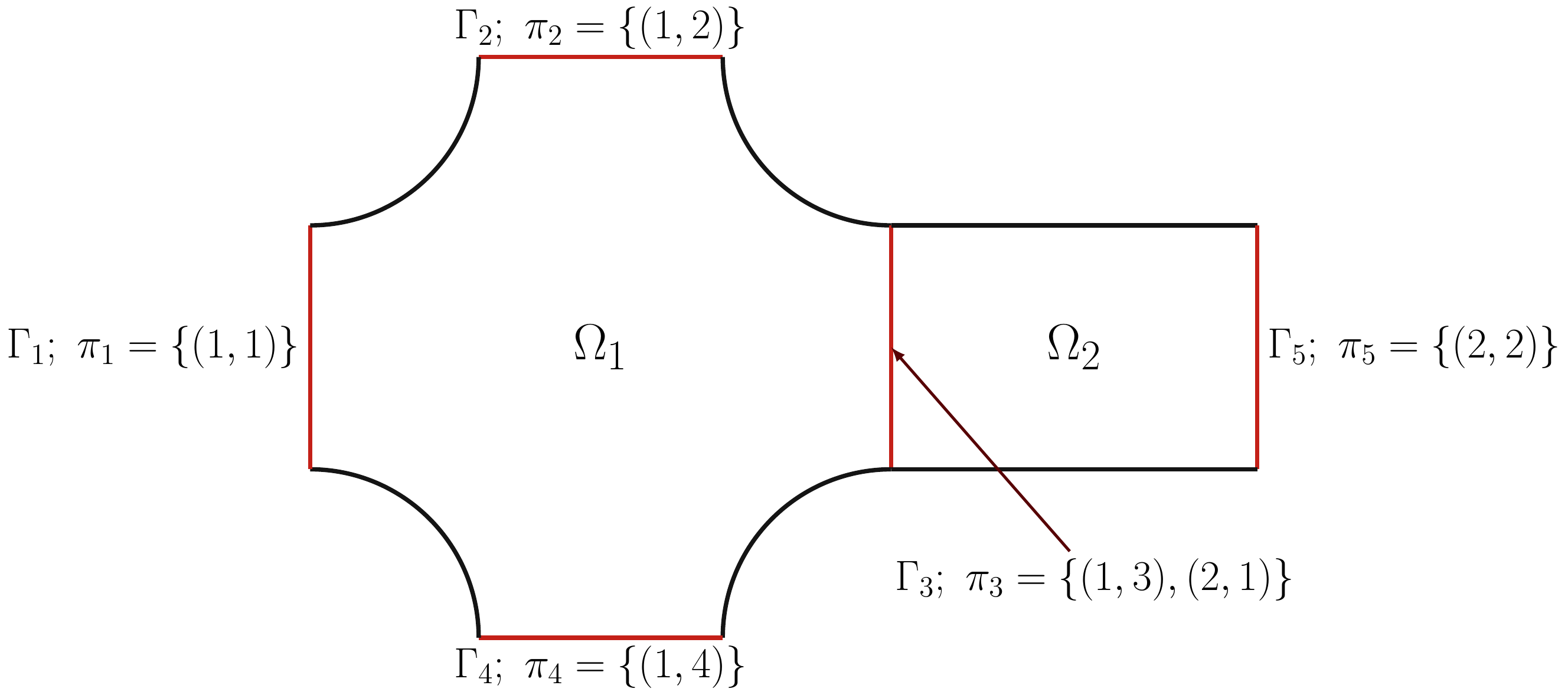}
    \caption{Illustration of PRSC nomenclature in a two-component system;
    red segments of the boundary are ports, with their global index and
    index set shown.}
    \label{fig:prsc_nomenclature}
\end{figure}

Now $\fulldomain$ is fully defined by specifying the component
domains $\domaini, \ i \in \{1, \ldots, \ncomps\}$, and the connections between them
in the form of port pairs $\left(\portkk{i}{j}, \portkk{i'}{j'}\right)$. Although
the members of this port pair refer to subsets of $\domaini$ and $\domaink{i'}$
(``local ports'') respectively, since they coincide in space we also refer to
them as a single ``global port''. A single, unconnected port also corresponds to a
single global port. Global ports are denoted by
$\globportp, \ p \in \{1,\ldots,\nglobports\}$ and defined by an index set
$\indexset = \left\{(i, j), (i', j')\right\}$, for connections of two local ports,
or $\indexset = \left\{(i, j)\right\}$ for disconnected ports. The correspondence
between global and local ports is illustrated in Fig.~\ref{fig:prsc_nomenclature}.
We also define a mapping $\portmapi$ such that given a local port $\portij$,
$\portmapi(j)$ maps $j$ to the index $p$ of the global port $\globportp$. 

The restriction of $X$ to $\portij$ (its ``port space'')
is denoted $\portspaceij$.  Where a global port $\globportp$
is the coincidence of two local ports,
$\indexset = \left\{(i, j), (i', j')\right\}$,
their port spaces are identical: $\portspaceij = \portspacekk{i'}{j'}$.
Letting the dimension of $\portspaceij$ be $\nportdofsij$, define a basis
$\left\{\pbasisijk \in \portspaceij,
\ k \in \left\{1, \ldots, \nportdofsij\right\}\right\}$
for $\portspaceij$, and a ``lifted port basis''
$\left\{\lpbasisijk \in \femspace|_{\domaini}, \ k \in
\left\{1,\ldots,\nportdofsij\right\}\right\}$, where
$\lpbasisijk|_{\portij} = \pbasisijk$ and
$\lpbasisijk|_{\portik{j'}} = 0, \ j' \neq j$. We lift the port bases by solving
\begin{equation}\label{eq:affine_lifting}
    \begin{split}
        \bilini
                \left(
            \lpbasisijk, v; \prscparamk{i;0}
        \right) &= 0,\ \forall v \in \bubspacei \\
        \lpbasisijk &= \pbasisijk \text{ on } \portij \\
        \lpbasisijk &= 0 \text { on } \portkk{i}{j'}, \ j' \neq j.
    \end{split}
\end{equation}

The global lifted port basis associated to
$\globportp$ is given by $\gpbasispk = \lpbasisijk + \lpbasiskkk{i'}{j'}$
if $\indexset = \left\{(i,j),(i',j')\right\}$, or $\gpbasispk = \lpbasisijk$ if
$\indexset = \left\{(i, j)\right\}$, where port basis functions are
extended by zero outside of their domain. 

In addition to port spaces, each component has a ``bubble space'' defined as
\begin{equation}
    \bubspacei \equiv
    \left\{
        \testfun \in \femspacei : \testfun|_{\portij} = 0,
        \ j \in \left\{1,\ldots,\nportsi\right\}
    \right\}.
\end{equation}
Define two classes of parameter dependent ``bubble functions'' by solution
of the equations
\begin{equation}\label{eq:port_space_bubbles}
    \bilini
    \left(
        \bubijk\left(\prscparami\right), v; \prscparami
    \right) =
    -\bilini
    \left(
        \lpbasisijk, v; \prscparami
    \right), \ \forall v \in \bubspacei
\end{equation}
\begin{equation}\label{eq:forcing_bubble}
    \bilini
    \left(
        \bubfi\left(\prscparami\right), v; \prscparami
    \right) =
    \lin_i\left(v; \prscparami\right), \ \forall v \in \bubspacei
\end{equation}
where $\bubijk, \bubfi \in \bubspacei$.
With these definitions, the solution to \eqref{eq:global_weak_form}
on $\domaini$ may be written as:
\begin{equation}\label{eq:cw_prsc_solution}
    \femsol(\prscparam)|_{\domaini} =
        \bubfi\left(\prscparami\right) +
        \sum_{j=1}^{\nportsi} \sum_{k=1}^{\nportdofsij}
            \prscsolkk{\portmapi(j)}{k}\left(\prscparam\right)
            \left(
                \bubijk(\prscparami) + \lpbasisijk
            \right)
\end{equation}
with $\prscsolkk{p}{k}$ unknown coefficients. 

Defining interface functions $\ifuncijk(\prscparami)$ by
\begin{equation}\label{eq:ifunc_def}
    \ifuncijk(\prscparami) \equiv \lpbasisijk + \bubijk(\prscparami)
\end{equation}
and global interface functions $\Ifuncpk(\prscparam)$ by
\begin{equation}\label{eq:glob_ifunc_def}
    \Ifuncpk(\prscparam) \equiv \sum_{(i, j) \in \indexset} \ifuncijk(\prscparami),
\end{equation}
the global solution to \eqref{eq:global_weak_form} may be written as
\begin{equation}\label{eq:global_prsc_solution}
    \femsol(\prscparam) =
        \sum_{i=1}^{\ncomps} \bubfi\left(\prscparami\right) +
        \sum_{p=1}^{\nglobports} \sum_{k=1}^{\nportdofsp}
        \prscsolkk{p}{k}\left(\prscparam\right) \Ifuncpk\left(\prscparam\right).
\end{equation}
Substituting \eqref{eq:global_prsc_solution} in \eqref{eq:global_weak_form} and
restricting test functions to the space spanned by the interface functions
leads to the condensed system of equations
\begin{equation}\label{eq:prsc_linear_system}
    \scstiff(\prscparam)\scsol(\prscparam) = \scforce(\prscparam),
\end{equation}
where the entries of $\scstiff$ and $\scforce$ are
\begin{equation}\label{eq:scstiffness_def}
    \scstiffk{(p,k), (p', k')} = \bilin
    \left(
        \Ifuncpk\left(\prscparam\right),
        \Ifunckk{p'}{k'}\left(\prscparam\right);
        \prscparam
    \right)
\end{equation}
and
\begin{equation}
    \scforcek{(p,k)} = \lin
    \left(
        \Ifuncpk\left(\prscparam\right); \prscparam
    \right) -
    \sum_{i=1}^{\ncomps} \bilin
    \left(
        \bubfi(\prscparami), \Ifuncpk\left(\prscparam\right);
        \prscparam
    \right).
\end{equation}
A pair $(p, k)$ identifies a single degree of freedom in
the linear system \eqref{eq:prsc_linear_system}:
the coefficient of $\Ifuncpk$ in \eqref{eq:global_prsc_solution}.
The Schur complement consists of local contributions from each component,
analagous to element stiffness matrices in a FEM. These contributions are
defined by
\begin{equation}\label{eq:component_scstiffness}
    \sccompstiffik{(j, k), (j', k')} = \bilini
    \left(
        \ifuncijk(\prscparami), \ifunckkk{i}{j'}{k'}(\prscparami);
        \prscparami
    \right).
\end{equation}
We make use of this definition in sensitivity analysis.

In this description we have referred only to a collection of $\ncomps$
different subdomains $\domaini$ and termed these ``components''. In practice a
distinction is made
between \textit{reference} and \textit{instantiated} components. Each instantiated
component is associated to a corresponding reference component, where the number
of reference components is small, but $\ncomps$ may be very large.
An offline data set needs only to be constructed for the small number of
reference components, and is used to construct the linear system \eqref{eq:prsc_linear_system}
using efficient operations on the offline data set.
To simplify the presentation of PRSC we have avoided introducing separate notation for
instantiated and reference components and describing the correspondence between the
two, but these details are key to the efficient implementation of the method.

This description of PRSC so far only includes static condensation.
Model reduction is incorporated by using a reduced port basis
$$\left\{\rpbasisijk \in \portspaceij, \ k \in \{1,\ldots,\nrportdofsij\}\right\},$$
with reduced port basis dimension $\nrportdofsij < \nportdofsij$,
and corresponding lifted port basis functions $\tilde{\psi}_{i,j,k}$. 
In presenting the Schur complement system, we have used the notation for
the full order model; however, the reduced order system may be derived by substituting the
reduced lifted port basis and its dimension wherever the full order lifted
port basis and its dimension appear. Henceforth, ROM quantities
are denoted by adding a tilde to their FOM counterpart; for example,
$\rfemsol$ is the ROM solution for displacement.
We determine an appropriate reduced port
basis using the pairwise training procedure described in \cite[Algorithm 2]{prsc}
and proper orthogonal decomposition with respect to the $L^2$ inner product
on each port.

\section{A component-wise formulation for topology optimization with stress constraints}
\label{sec:cw_formulation}
Here we present the main contribution of this work: a method for stress-based TO that
leverages PRSC to enable an efficient optimization, even when using
a high-resolution discretization. With the use of PRSC, we also inherit its
greatest advantage, which is the ability to reuse the component library constructed
in the offline phase to model any arrangement of instantiated components.
This allows the use of the same offline dataset to solve many TO problems in
various geometries.

Below, we describe our methodology in detail: the formulation of the
optimization problem, the material model and parameterization, the
form of the aggregated stress constraints, and finally, postprocessing
considerations.

\subsection{Optimization formulation}
Our goal is to solve the problem
\begin{equation}
\begin{aligned}\label{eq:basic_formulation}
  & \underset{\densityvec \in \{0,1\}^{\ncomps}}{\text{minimize}} \hspace{5pt} \objective(\densityvec) =
  \sum_{i = 1}^{\ncomps} \densityi \left|\domaini\right| \\
  & \text{s.t.} \max_{\Omega} \hspace{5pt} \stress(\densityvec) \leq \stressmax,
\end{aligned}
\end{equation}
where $\densityvec \in \{0, 1\}^{\ncomps}$ is the parameter vector, containing a density
for each component in the system with $\densityi = 0$ meaning that component $i$ is removed
and $\densityi = 1$ meaning that component $i$ is fully solid;
$\left|\domaini\right|$ is the volume of $\domaini$; $\objective(\densityvec)$ gives
the mass of the system; and $\stress$ is some yield criterion.
In this work, $\stress$ is the Von Mises yield stress. $\stressmax$ is an
upper bound on the allowed value of the yield criterion $\stress$.

To render \eqref{eq:basic_formulation} solvable by a gradient-based optimization
solver, we make two changes following common practice.
The first is to allow the component densities
to reside in a continuum: $\densityvec \in [0, 1]^{\ncomps}$. Material properties
for intermediate densities are interpolated by a standard SIMP scheme.
Second, because of the non-differentiable nature of the max function, we replace the
constraint in \eqref{eq:basic_formulation} with a differentiable approximation.
With these modifications, the optimization formulation becomes:
\begin{equation}
  \begin{aligned}\label{eq:actual_formulation}
    & \underset{\densityvec \in [0, 1]^{\ncomps}}{\text{minimize}} \hspace{5pt} \objective
    (\densityvec) = \sum_{i=1}^{\ncomps} \densityi \left|\domaini\right| \\
    & \text{s.t.}\ \constrainti(\densityvec) \leq 1, \ 1 \leq m \leq \naggs,
  \end{aligned}
\end{equation}
with the aggregate constraints $\constrainti$ defined in Section \ref{sec:stress_agg}.
This is the optimization formulation that we actually solve. A local minimum
of \eqref{eq:actual_formulation} provides an approximation to solution of
\eqref{eq:basic_formulation}, but the challenging integer programming problem
has been replaced by one solvable using gradient-based
nonlinear optimization solvers.

\subsection{Material model \& component-based parameterization}
We model a material described by isotropic linear elasticity.
On each component, the bilinear form $\bilini$ is given by 
\begin{equation}\label{eq:linel_bilinear_form}
    \bilini(\femsol, \testfun; \densityi) = \int_{\domaini}
    \simpscaling(\densityi)
    \elastictensor(\femsol) \cdot \nabla \testfun \ dx,
\end{equation}
where $\simpscaling$ is the SIMP penalization term and $\elastictensor$
is the symmetric elasticity tensor given by
\begin{equation}\label{eq:elastic_tensor}
    \elastictensor(\femsol) = \lamemu
    \left(\nabla\femsol + \nabla\femsol^T\right)
    + \lamelambda(\nabla \cdot \femsol)\Identity,
\end{equation}
with $\lamelambda$ and $\lamemu$ the Lam\'{e} parameters. Note that $\lamemu$
is a fixed material parameter and
bears no relation to the vector of component parameters $\prscparami$ from the
discussion of PRSC. In our optimization formulation
$\prscparami = \left[\densityi\right]$, so we use the scalar parameter $\densityi$
in its place.
This form for $\bilini$ satisfies the
necessary assumptions to implement the affine simplification of PRSC described
in our previous work \cite{our_first_paper}. This simplification allows us to
compute Schur complement contributions offline, making the model even more
economical than PRSC without this simplification.

The linear form $\lin_i$ is
\begin{equation}\label{eq:linel_linear_form}
    \lin_i(\testfun) = \int_{\domaini} \force \cdot \testfun\ dx,
\end{equation}
where $\force$ is the vector of forces applied per unit volume. We
assume $\force$ to be independent of parameter in this work; this assumption
is not necessary, but simplifies the sensitivity analysis.

We aim to obtain a ``black and white'' solution, $\densityi \in \left\{0, 1\right\}$.
To drive the continuous density values to black and white
solutions, we use the solid isotropic material (SIMP) parameterization
\cite{bendsoe_1989,zhou_rozvany_1991}. The Young's modulus is
penalized using a power law in the density: $\youngsi \propto \densityi^p$,
with $\youngsi$ the Young's modulus of the material in $\domaini$.
In our optimizations, we choose $p = 3$. To ensure a well-posed problem,
the SIMP parameterization is modified so
that $s(\densityi)$ in \eqref{eq:linel_bilinear_form} is given by
\begin{equation}\label{eq:simp}
    s(\densityi) = \left[\densityi + (1 - \densityi) \densitymin\right]
    ^3,
\end{equation}
where $\densitymin > 0$ is chosen to be small ($10^{-3}$). This ensures that the
linear system \eqref{eq:prsc_linear_system} possesses a solution,
but the stiffness of components with $\densityi = 0$ is negligible compared
to components with $\densityi = 1$.

This parameterization is the same as a typical topology optimization
using element-wise density variables, but enables the use of a component-wise
ROM to reduce the computational cost because each component ROM only depends
on a single density parameter. In a conventional, element-based density TO,
the large dimension of parameter space makes model reduction impractical.
The component-wise scheme has an additional benefit: it obviates the need
for a density filter (as in
\cite{le_stressbased_2010,holmberg_stressbased_2013,senhora_et_al_2020,luo_kang_2012}
and others) to avoid checkerboarding. The imposition of a ground structure
intrinsically imposes a length scale without the use of a filter.

\subsection{Imposing stress constraints using the ROM}
Our formulation of stress constraints consists in three parts:
stress relaxation, aggregation using the KS functional, and rewriting
of the aggregate constraints using ROM operators. Of these parts, the first
two are standard practice; the rewriting using ROM operators is an original
contribution.

\subsubsection{Stress relaxation}
It is well known that stress constraints in the absence of some relaxation
result in an optimization problem where minima are contained in a
degenerate subspace (the ``singularity problem'') \cite{kirsch_1990}.
As remarked in \cite{duysinx_sigmund_1998}, this occurs where material
in some region must vanish to reach a local optimum, but the material
remains strained so that solutions for intermediate densities in that
region violate constraints. When material has vanished, however, these
constraints should be ignored -- the constraint is discontinuous. Stress
relaxation smooths this discontinuity to allow a gradient-based optimization
algorithm to reach solutions in the degenerate subspace.

We adopt the $qp$-relaxation proposed by Bruggi \cite{bruggi_2008}, with
the particular choice of $q = 2.5$ as in
\cite{le_stressbased_2010,holmberg_stressbased_2013}. Denoting the Von Mises
yield stress by $\vonmises$, the relaxed stress is then given by
\begin{equation}
    \relaxedstress|_{\domaini} = s'(\densityi) \vonmises,
\end{equation}
where $\vonmises$ is the Von Mises stress as computed using the stiffness
of the base material (not the SIMP penalized stiffness), and 
\begin{equation}\label{eq:relaxation}
    s'(\densityi) = \densityi^{1/2}.
\end{equation}
We note that the optimization has the trivial solution
$\densityvec = \boldsymbol{0}$; in practice, this solution is not reached by
the optimizer. Convergence to the trivial solution, if problematic, could be
addressed by imposing a lower bound on $s'(\densityi)$ in the same fashion
as in \eqref{eq:simp}. Applying stress relaxation, the stress constraint
becomes
\begin{equation}\label{eq:relaxed_max_constraint}
    \max_{\Omega}\ \relaxedstress(\densityvec) \leq \stressmax,
\end{equation}
which is identical to the constraint in \eqref{eq:basic_formulation} for
black and white solutions.

\subsubsection{Stress aggregation}\label{sec:stress_agg}
We cannot use the constraint \eqref{eq:relaxed_max_constraint} in a
gradient-based optimization because the max function is non-differentiable.
To circumvent this difficulty, we use stress aggregation to provide a
differentiable approximation of the max function. Aggregation strategies include the
$p$-norm, the Kreisselmeier-Steinhauser (KS) functional, or an ``induced aggregation''
functional as described by Kennedy and Hicken \cite{kennedy_hicken_induced_2015}. We
use a continuous KS aggregation, with stresses aggregated over multiple
aggregation domains (as  in
\cite{holmberg_stressbased_2013,le_stressbased_2010}, and others).
For numerical stability
in finite precision arithmetic, we aggregate the ratio of the relaxed
stress to the maximum stress, rather than the relaxed stress itself.
The single constraint \eqref{eq:relaxed_max_constraint} becomes $\naggs$ constraints
given by
\begin{equation}\label{eq:agg_constraints}
    \constrainti(\densityvec) = \frac{1}{p} \ln
    \left(
        \frac{1}{\alpha}\int_{\aggdomaini} \exp
        \left(
            p \frac{\relaxedstress(\densityvec)}{\stressmax}
        \right) dx
    \right) \leq 1, \ 1 \leq m \leq \naggs, 
\end{equation}
where $p$ and $\alpha$ are fixed parameters of the aggregation. Increasing $p$ results in a closer approximation of
the max function, but also a more difficult optimization problem.
$\alpha$ is a normalization, whose determination we address below. 
The aggregation domains $\aggdomaini$ are not the same as component domains
$\domaini$, nor even spatially contiguous domains.
Their determination is addressed in the next section.
Ideal values of $p$ and $\naggs$ are problem dependent, and studied in our numerical
experiments.

The normalization $\alpha$ is
computed based on the following observation from \cite[][Eq. 8]{kennedy_hicken_induced_2015}:
for given values of $\densityvec$ and $p$, we will have
\begin{equation*}
    \constrainti(\densityvec) \geq \max_{\aggdomaini} \frac{\relaxedstress(\densityvec)}{\stressmax}
\end{equation*}
if
\begin{equation}\label{eq:alpha_estimation}
    \alpha \leq \int_{\aggdomaini} \exp\left(
        \frac{p}{\stressmax}\left[\relaxedstress - \max_{\aggdomaini}\relaxedstress\right]
    \right)dx.
\end{equation}

Let the initial value of $\densityvec$ in the optimization be $\densityvec_0$.
Using \eqref{eq:alpha_estimation}, we compute values of $\alpha$ for each aggregation
domain at the initial condition of the optimization as follows: compute
\begin{equation}\label{eq:alpha-determination}
    \alpha_k = 
    \exp\left(
        \frac{p}{\stressmax}\left[
          \max_{\aggdomaini}\relaxedstress(\densityvec_0)\right]
        \right)^{-1}
    \int_{\aggdomaini} \exp\left(
        \frac{p}{\stressmax}\relaxedstress(\densityvec_0)
    \right) dx,
\end{equation}
then set $\alpha$ for all aggregates to be equal to the minimum of all $\alpha_k$.

\subsubsection{Achieving a conservative optimization result}\label{sec:conservative_result}
The aggregation described in the previous section is not necessarily conservative;
for some values of $\densityvec$, the
constraints \eqref{eq:agg_constraints} may be satisfied while the true
max stress constraint \eqref{eq:relaxed_max_constraint} is not. To consider an optimization
successful, we require satisfaction of the latter. In most cases determining $\alpha$ from Eq. \eqref{eq:alpha-determination} is
sufficient to obtain a conservative optimization solution; where this is not the case,
we use a simple heuristic:
substitute $\heurstressmax < \stressmax$ for $\stressmax$ in
\eqref{eq:agg_constraints}. The choice of $\heurstressmax$ is problem dependent.
Values used are reported in the numerical experiments -- for most,
$\heurstressmax = \stressmax$.

This heuristic approach has mathematical justification. Eq. \eqref{eq:agg_constraints}
may be rewritten as
\begin{equation}
    \constrainti(\densityvec) = \frac{1}{p}\ln\left[
        \int_{\aggdomaini} \exp\left(
            p\frac{\relaxedstress(\densityvec)}{\stressmax}
        \right)dx \right]- \frac{\ln \alpha}{p} - 1 \leq 0.
\end{equation}
Thus if $\alpha < e$, the value of the constraint is increased by a constant factor.
A decrease in the upper bound in the optimization likewise corresponds to a constant
increase in the value of $\constrainti$; therefore the imposition of $\heurstressmax
< \stressmax$ may be viewed as a correction to the value of $\alpha$ estimated from
\eqref{eq:alpha_estimation}, necessary because that value is conservative only for
a particular value of $\densityvec$.

We also assign aggregation regions differently from previous works. The two main
strategies are to either distribute stress evenly among aggregation regions
\cite{le_stressbased_2010}, or to assign the highest
stresses to a single aggregation region, resulting in a better approximation of the
maximum in that region \cite{holmberg_stressbased_2013}.
However, because we do not adaptively reassign regions and cannot know
how stresses will be distributed at a given optimization iteration,
neither of these approaches is possible here.
Instead, we approximate the effect of an approach that computes an even
distribution of stresses by assigning aggregation regions randomly.
For each $\domaini$, we randomly assign an equal number of elements to each
of the $\naggs$ aggregation domains $\aggdomaini$ (illustrated in Fig.~\ref{fig:agg-regions}).
Ideally, this assignment of aggregation regions results in a problem where
stresses are evenly distributed between the different aggregation regions.
We show results on the same problem for different assignments of
aggregation regions in Section \ref{sec:agg_region_comparison}.

\begin{figure}
    \centering
    \includegraphics[width=0.5\textwidth]{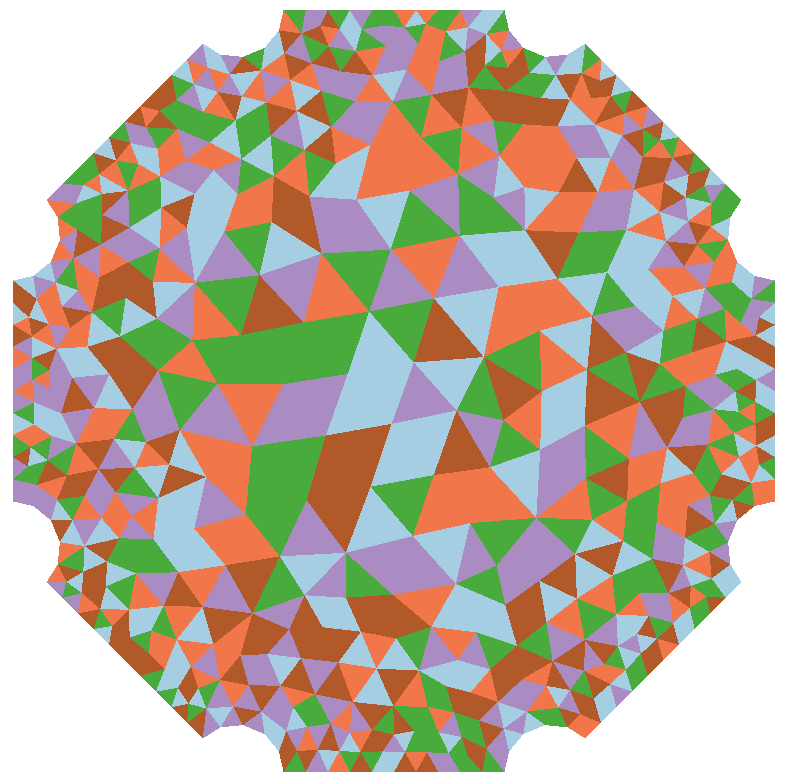}
    \caption{Illustration of aggregation region assignment for a single component and
    five aggregation regions. Elements contained in the same aggregation region are
    assigned the same color.}
    \label{fig:agg-regions}
\end{figure}

\subsubsection{Efficient computation of aggregates using the CWROM}\label{sec:fast_aggregation}
The aggregation in Eq.~\eqref{eq:agg_constraints} is nonlinear in $\vonmises$,
which is in turn a nonlinear function of displacement. This nonlinearity prohibits
the expression of the constraints in terms of $\rscsol(\densityvec)$ without
reconstructing the displacement field $\rfemsol(\densityvec)$.
However, computing the constraints by reconstructing $\rfemsol(\densityvec)$
and performing quadrature in a loop over elements is inefficient because of
the high dimension of the underlying discretization.
For example, in our smaller
numerical example of Section \ref{sec:l_bracket}, the mesh underlying the
component-wise ROM has about 2.5 million elements. Integration over each element
requires retrieving coefficients of the finite element basis functions, which
are generally not contiguous in memory, and computing values of the non-linear
integrand as a function of these coefficients. While it is not possible to remove
the dependence of the cost of aggregation on the FOM dimension, we reduce the
cost by expressing aggregation in terms of more cache-friendly operations using
$\rscsol(\densityvec)$ instead of $\rfemsol(\densityvec)$.

We rewrite the constraints in terms of the ROM by noting that $\vonmises$, which
is not linear in displacement, may be written as a function of the Cauchy stress tensor
components. Limiting ourselves to two dimensions, we define vectors
\begin{equation}\label{eq:stress-tensor-components}
    \stressvecixx\left(\rscsol(\densityvec)\right) = \stressopixx\rscsol,
    \ \stressveciyy\left(\rscsol(\densityvec)\right) = \stressopiyy\rscsol,
    \ \stressvecixy\left(\rscsol(\densityvec)\right) = \stressopixy\rscsol
\end{equation}
containing the values of stress tensor components at quadrature points
for aggregation region $\aggdomaini$. We assume that there are no forces applied
except on ports. Then, from the definition of the Cauchy stress tensor and
Eq.~\eqref{eq:global_prsc_solution}, the linear operators $\stressopixx$,
$\stressopiyy$, and $\stressopixy$ are given by:
\begin{align}
    \label{eq:stress-operators}
    \begin{split}
        \stressopixx &= (2\mu+\lambda) \ifuncxpx + \lambda \ifuncypy \\
        \stressopiyy &= (2\mu+\lambda) \ifuncypy + \lambda \ifuncxpx \\
        \stressopixy &= \mu\left(\ifuncxpy + \ifuncypx\right),
    \end{split}
\end{align}
with $\ifuncxpx$ a matrix such that $\left[\ifuncxpx\right]_{j,k}$ contains the partial derivative
with respect to $x$ of the $x$ component of the $k$-th interface function at
the $j$-th quadrature point in $\aggdomaini$, and $\ifuncypy$,
$\ifuncxpy$, and $\ifuncypx$ defined similarly.

In terms of the stress tensor component vectors \eqref{eq:stress-tensor-components},
the vector of Von Mises stresses at quadrature points in $\aggdomaini$ is given
by
\begin{equation}
    \vmquadi\left(\rscsol(\densityvec)\right) = \left(
        \hadamard{\stressvecixx}{\stressvecixx} +
        \hadamard{\stressveciyy}{\stressveciyy} -
        \hadamard{\stressvecixx}{\stressveciyy} +
        3\hadamard{\stressvecixy}{\stressvecixy}
    \right)^{1/2},
\end{equation}
where $\odot$ denotes the Hadamard product and the exponentiation is applied
element-wise. We omit the functional dependence of $\stressvecixx$, etc.
on $\rscsol(\densityvec)$ for brevity. The vector of relaxed stresses at
quadrature points is then
\begin{equation}
    \rstressveci\left(\densityvec, \rscsol(\densityvec)\right) = \hadamard{s'\left(\aggrhoi\right)}{\vmquadi},
\end{equation}
where $\aggrhoi$ contains the value of the density parameter at each quadrature
point and $s'$ is the relaxed stress \eqref{eq:relaxation}, applied
elementwise.

Taking $\boldsymbol{w}$ to be a vector containing the coefficients for
quadrature over $\aggdomaini$, we can rewrite the constraints in
Eq.~\eqref{eq:agg_constraints}:
\begin{equation}
    \label{eq:efficientks}
    \constrainti\left(\densityvec, \rscsol(\densityvec)\right) =
    \frac{1}{p}\ln\left(
        \frac{1}{\alpha}\boldsymbol{w}^T \left[\exp\left(
            \frac{p}{\stressmax} \rstressveci\left(\densityvec,
            \rscsol(\densityvec)\right)
        \right)\right]
    \right) - 1 \leq 0,
\end{equation}
with exponentiation interpreted elementwise.

Although this expression improves the performance of aggregation by
rewriting it in a more cache-friendly manner, the size of the stress operators
in Eq.~\eqref{eq:stress-operators} grows with the number of quadrature points
and the dimension of the ROM, until at some point it is more efficient to
compute the aggregates in the obvious fashion. In our numerical examples,
we use a four-point quadrature rule exact for second-order polynomials.
Also, because we assign aggregation regions randomly, the stress operators
in \eqref{eq:stress-operators} cannot be computed offline; instead they are
computed once at the beginning of the optimization. The cost of this computation
is insignificant relative to the total cost of optimization.

\subsubsection{Sensitivity analysis of stress constraints}\label{sec:sensitivity}
To derive the gradient of the constraint $\constrainti$ with respect to $\densityvec$,
we note that $\constrainti(\densityvec)$ depends both directly on
$\densityvec$ and implicitly through $\rscsol(\densityvec)$, with the latter dependence
defined by the forward model. Therefore, the gradient is given by a total derivative:
\begin{equation}
    \label{eq:total_derivative}
    \nabla_{\densityvec}\ \constrainti\left(
        \densityvec, \rscsol(\densityvec)
    \right) =
    \pderiv{\constrainti}{\densityvec} + 
    \left(\deriv{\rscsol}{\densityvec}\right)^T \pderiv{\constrainti}{\rscsol},
\end{equation}
where the notation $\deriv{\rscsol}{\densityvec}$ indicates the Jacobian matrix
with
\begin{equation*}
    \left[
        \deriv{\rscsol}{\densityvec}
    \right]_{i,j} =
    \pderiv{\rscsol_i}{\densityk{j}}.
\end{equation*}
The action of $\left(\deriv{\rscsol}{\densityvec}\right)^T$ is computed using the
adjoint method. From Eq. \eqref{eq:prsc_linear_system}, we have
\begin{equation}
    \deriv{\rscstiff}{\densityvec} \rscsol + \rscstiff \deriv{\rscsol}{\densityvec}
    =
    \deriv{\rscforce}{\densityvec}.
\end{equation}
In our problem settings the forcing is independent of parameter and the
above yields
\begin{equation}
    \deriv{\rscsol}{\densityvec} = -\rscstiff^{-1} \left[
        \deriv{\rscstiff}{\densityvec} \rscsol
    \right].
\end{equation}
By substitution in Eq. \eqref{eq:total_derivative}, we obtain the gradient:
\begin{equation}
    \label{eq:adjoint_sensitivity}
    \nabla_{\densityvec} \ \constrainti\left(
        \densityvec, \rscsol(\densityvec)
    \right)
    = \pderiv{\constrainti}{\densityvec} - \left[
        \deriv{\rscstiff}{\densityvec} \rscsol
    \right]^T \adjointsol,
\end{equation}
where, making use of the symmetry of $\rscstiff(\densityvec)$,
$\adjointsol$ is given by solution of the adjoint equation
\begin{equation}
    \label{eq:adjoint_solve}
    \rscstiff(\densityvec) \adjointsol(\densityvec) = \pderiv{\constrainti}{\rscsol}(\densityvec).
\end{equation}
Note that $\deriv{\rscstiff}{\densityvec}$ is a third-order tensor; for
this sensitivity analysis, its product with $\rscsol$ is defined by:
\begin{equation}
    \left[
        \deriv{\rscstiff}{\densityvec} \rscsol
    \right]_{ij} =
    \sum_k \deriv{\rscstiff_{ik}}{\densityk{j}} \rscsol_k
\end{equation}

To close Eq. \eqref{eq:total_derivative}, we require expressions for the partial derivatives
of $\constrainti$. These are found using the chain rule: 
\begin{align}
    \pderiv{\constrainti}{\densityvec} &= \left(
        \pderiv{\rstressveci}{\densityvec}\right)^T
        \pderiv{\constrainti}{\rstressveci}
    \label{eq:dgdrho} \\
    \pderiv{\constrainti}{\rscsol} &= \left(
        \pderiv{\rstressveci}{\rscsol} \right)^T
        \pderiv{\constrainti}{\rstressveci}.
        \label{eq:dgdU}
\end{align}
Defining the vector
\begin{equation}
    \boldsymbol{e}\left(
        \densityi, \rscsol(\densityvec)
    \right) =
    \exp\left[
        \frac{p}{\stressmax}\rstressveci\left(
            \densityi, \rscsol(\densityvec)
        \right)
    \right],
\end{equation}
we find from Eq. \eqref{eq:efficientks} that 
$\pderiv{\constrainti}{\rstressveci}$ is
\begin{equation}
    \pderiv{\constrainti}{\rstressveci} = 
    \frac{1}{\stressmax \ \boldsymbol{w}^T
    \boldsymbol{e}\left(\densityvec, \rscsol(\densityvec)\right)}
    \ \hadamard{
        \boldsymbol{e}\left(
            \densityvec, \rscsol(\densityvec)
        \right)}{\boldsymbol{w}}.
\end{equation}

Finally, we require the derivatives of
$\rstressveci$ in Eqs. \eqref{eq:dgdrho} and \eqref{eq:dgdU}. With respect to
$\densityvec$, obtain
\begin{equation}
    \pderiv{\rstressveci}{\densityvec} =
    \prolongi\left(\hadamard{\deriv{s'}{\rho}(\aggrhoi)}{\vmquadi}\right) \otimes I,
\end{equation}
with $\prolongi$ a prolongation operator taking values in $\aggrhoi$ to their
corresponding entry in $\densityvec$ and padding with zeros for entries in
$\densityvec$ not corresponding to elements in $\aggdomaini$, and $\otimes$
denoting a row-wise Kronecker product.

With respect to $\rscsol$, the Jacobian of $\rstressveci$ is given by
\begin{equation}
    \pderiv{\rstressveci}{\rscsol} = s'\left(\aggrhoi\right) \otimes
    \pderiv{\vmquadi}{\rscsol},
\end{equation}
with $\pderiv{\vmquadi}{\rscsol}$ given by
\begin{equation}
    \pderiv{\vmquadi}{\rscsol} =
    \frac{1}{2} \left(\vmquadi\right)^{-1} \otimes \left[
        2 \stressvecixx \otimes \stressopixx +
        2 \stressveciyy \otimes \stressopiyy -
        \stressvecixx \otimes \stressopiyy -
        \stressveciyy \otimes \stressopixx +
        6 \stressvecixy \otimes \stressopixy
    \right].
\end{equation}

\subsection{Postprocessing}
\label{sec:postprocessing}
SIMP penalization does not necessarily succeed at creating black and
white solutions. Therefore, we postprocess by removing components with
$\densityi$ less than a prescribed minimum value $\densitymin$ from the
domain entirely, and setting the densities of remaining components to unity.
This reveals some amount of wasted mass in the structure where a component
that was previously attached to another now has a ``hanging'' port not on
a load-bearing path. In a post-processing step, we substitute streamlined
versions of these components from a library of alternatives, as illustrated
by Fig.~\ref{fig:component-deletion}. This procedure decreases the mass of
the system from that at the local optimum, often substantially, and its
implementation is simple in the component-wise framework.

\begin{figure}
    \centering
    \includegraphics[width=0.5\textwidth]{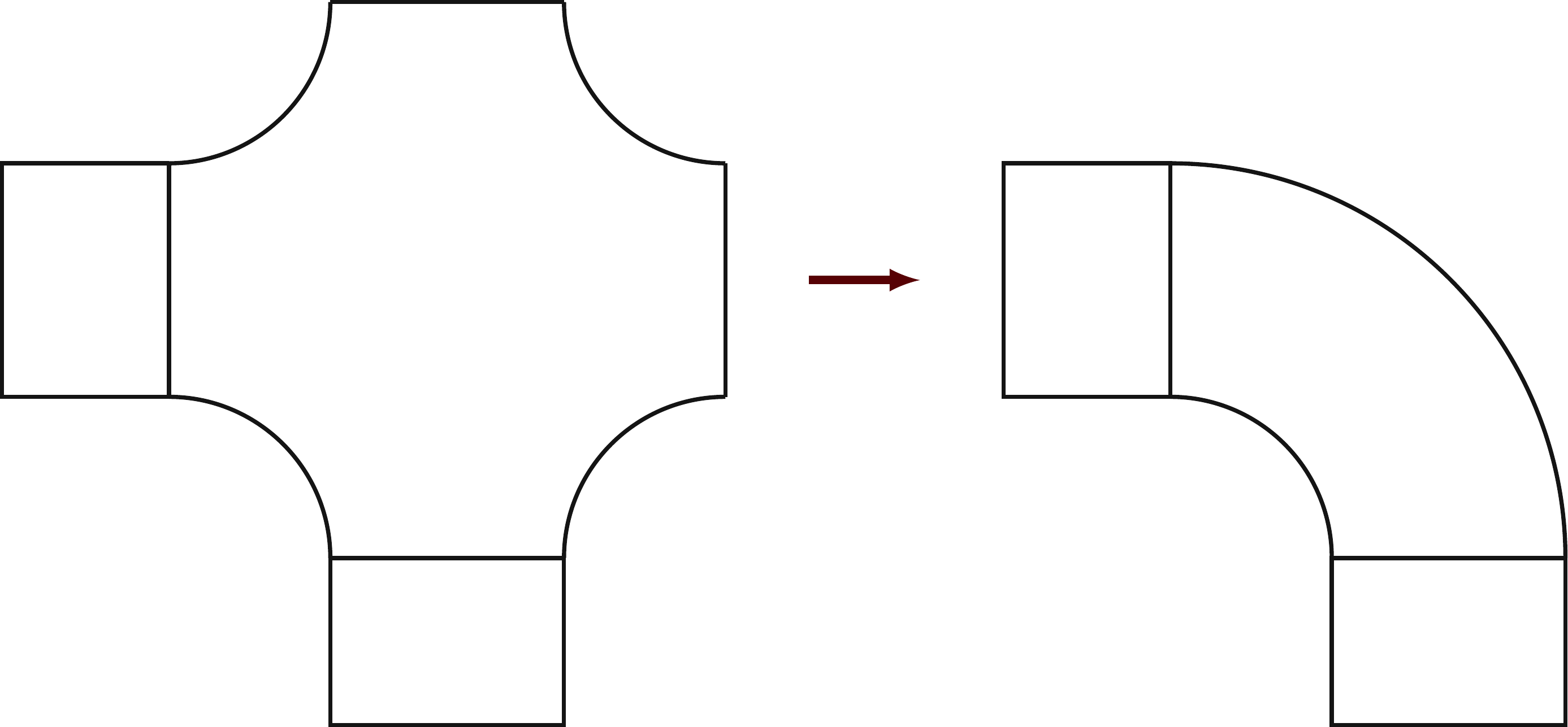}
    \caption{Illustration of replacing a component with two hanging ports by
             an alternative, streamlined geometry}
    \label{fig:component-deletion}
\end{figure}

\section{Numerical results}\label{sec:numerical_results}
We present numerical results for a pair of two-dimensional stress constrained problems.
In the first, we minimize mass
of an L-shaped bracket subject to a vertical load on its top right tip; in the second,
we minimize the mass of a cantilever beam, also with a vertical load at the top
right. The examples share the same material
properties: Young's modulus of 113.8 GPa, Poisson's ratio of 0.34, and a Von Mises yield criterion
of $\stressmax = 880$ MPa. Linear elasticity is approximated in two dimensions assuming plane stress conditions.
We demonstrate that our methodology succeeds in creating designs
that respect a constraint on maximum stress while approximating that constraint
using PRSC and stress aggregation.

\subsection{Problem setup}

\begin{figure}
    \centering
    \includegraphics[width=0.85\textwidth]{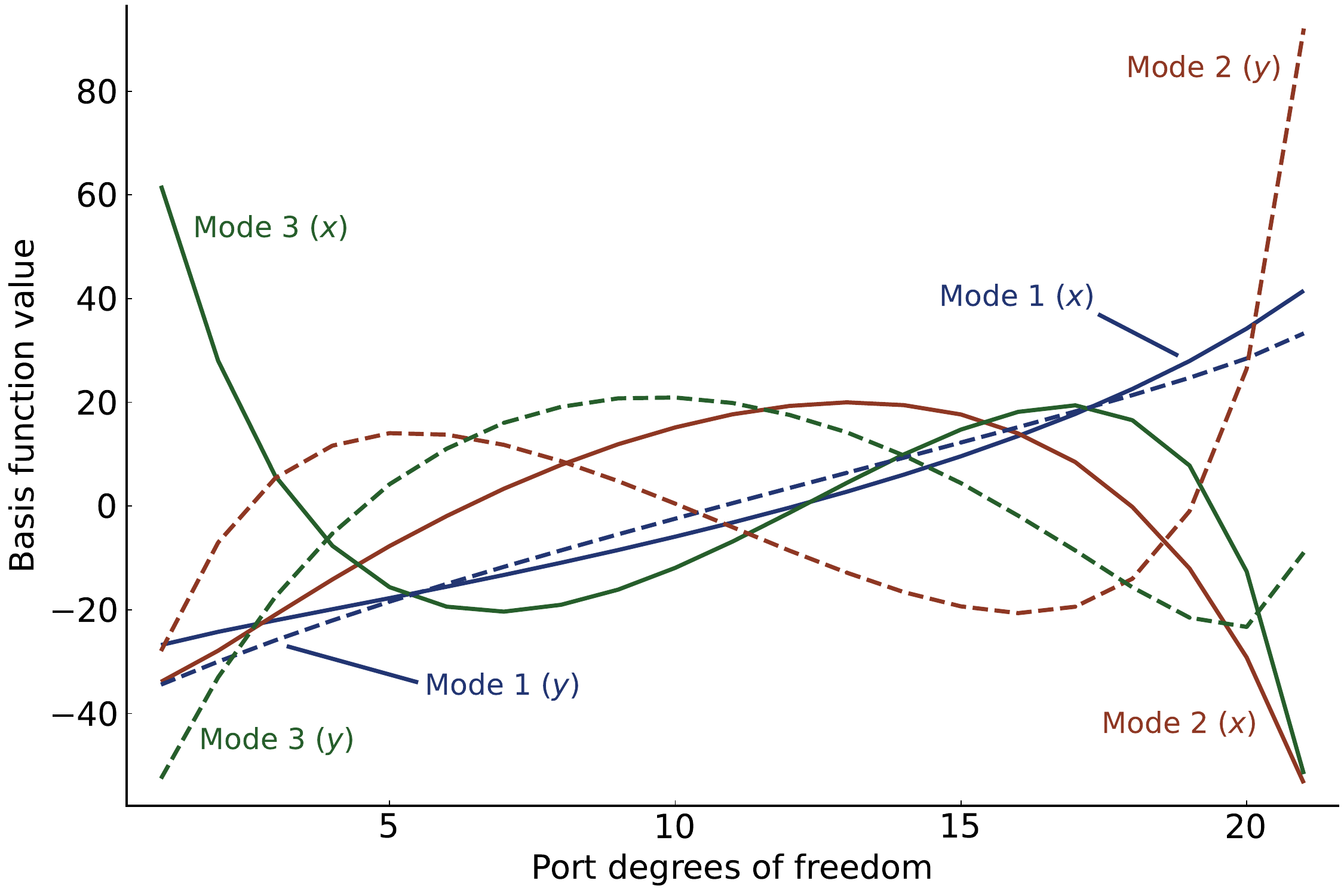}
    \caption{An example port basis, not including the constant function.
    The basis functions for the $x$ component of displacement are shown using
    a solid line; for the $y$ component, a dashed line.}
    \label{fig:port-basis-example}
\end{figure}

The reduced order model used for our primary optimization results is trained by collecting
500 snapshots for each pair of connected reference ports and constructing port
bases that capture 99.9\% of the total energy of each snapshot set, unless a different
total energy percentage is specified. When training,
we compute the POD for the X and Y components of displacement separately but using the
same snapshot set, and apply the total energy criterion to the POD for each component
of displacement separately to determine bases.
The regularization parameter for pairwise training was set to $\eta = 10$ \cite[Sec. 3.2.2]{prsc}.
The resulting bases contain four basis functions, including the constant function,
for each component of displacement. An example reduced basis is pictured in
Fig.~\ref{fig:port-basis-example}. For our primary optimization results, the
ROM is constructed using second-order triangular
finite elements; the full basis for a port contains twenty-one basis functions for each
component of displacement. We address the impact of using second-order elements vs.
first-order in Section \ref{sec:first_order_comparison}.

All optimizations are initialized from a fully solid ground structure, i.e., $\densityvec =
\boldsymbol{1}$. When postprocessing, we validate the optimized design by performing
analysis with the full-order model. Maximum stress values in postprocessing are
computed using values at cell centers, and at the quadrature points of a nine-point
quadrature rule in each element. We do not consider values of stress on
the boundary of an element due to its discontinuity there.

The optimization problem is solved using an interior point method as
implemented in Ipopt \cite{ipopt}. We found it was necessary to set the
parameter \verb|theta_max_fact| to a small value (0.5 was used in these results)
to prevent large steps into infeasible regions from which it is difficult
for the optimizer to escape. For the same reason, we disable the watchdog
procedure that attempts to escape regions of slow convergence by disabling
line search for one iteration.
In all of the optimization runs presented, Ipopt's convergence tolerance
was set to $10^{-6}$.

The Schur complement system is solved using a Cholesky factorization
as implemented in Eigen \cite{eigenweb}, and all optimizations are performed
without use of parallelism, for reproducibility. Timings are performed on a Linux
server with 2 AMD EPYC 7H12 processors with 64 physical cores and a base clock speed
of 2.6 GHz, and 2 TB of memory.

\subsection{Mass minimization of an L-bracket}
\label{sec:l_bracket}

Our first example addresses the design of an L-shaped bracket, a standard benchmark
in stress-based TO. While in most TO approaches the L-bracket presents a challenge because
of its reentrant corner, using the component-wise methodology,
we eliminate this difficulty by a careful choice of the ground structure.

\begin{figure}
    \centering
    \includegraphics[width=0.75\textwidth]{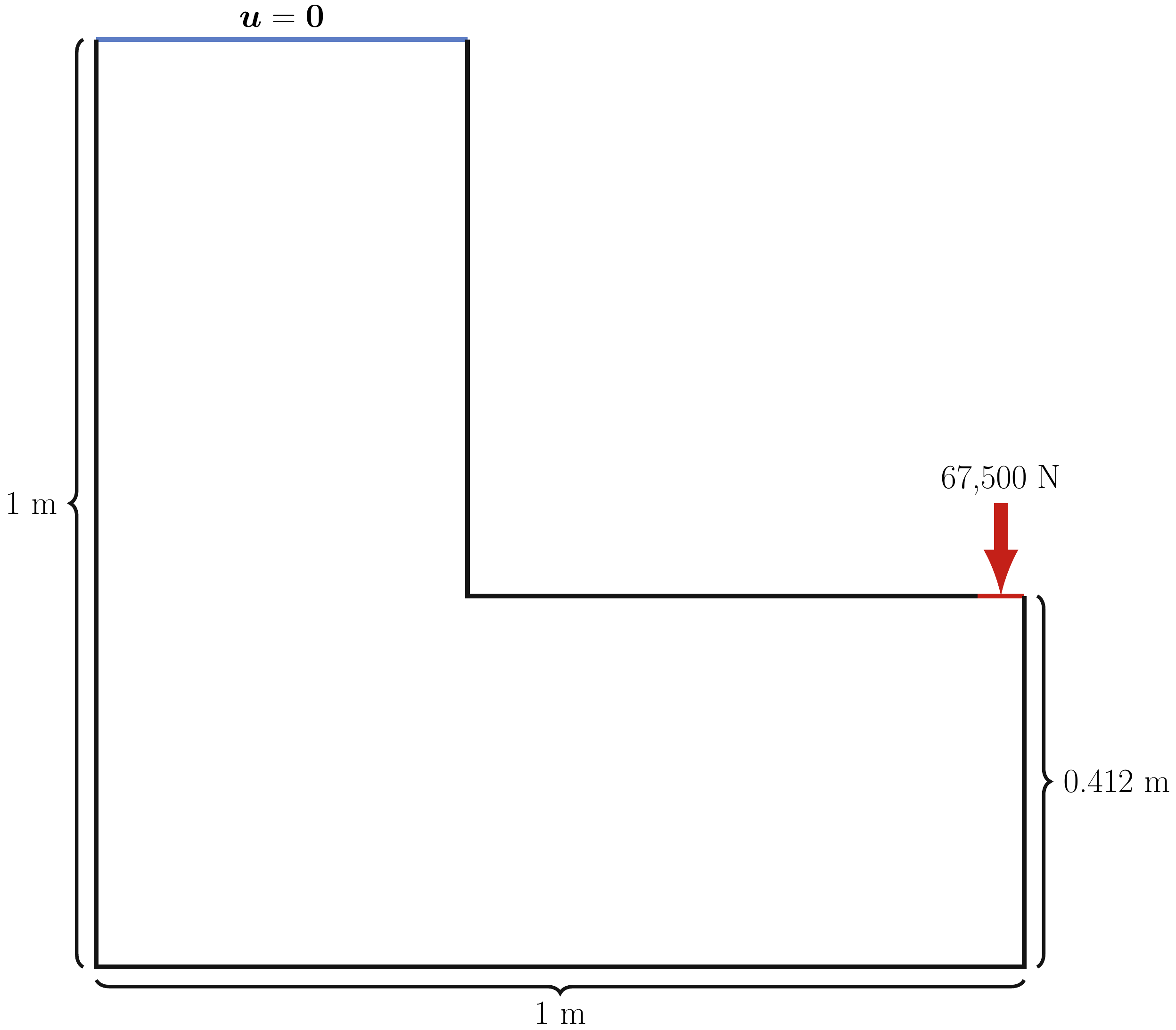}
    \caption{Illustration of the L-bracket geometry and loading condition.}
    \label{fig:l-bracket-setup}
\end{figure}

The problem setup for the L-bracket optimization is illustrated in Figure~\ref{fig:l-bracket-setup}. The thickness of the bracket into the page is
taken to be 5 cm. The ground structure contains 4,257 components and
2,427,362 finite elements; the reduced static condensation system
contains 44,214 unknowns while the full order static
condensation system has 278,544. The underlying finite element model contains more
than 5 million unknowns (not accounting for those degrees of freedom on ports, which
are counted twice). The bracket's upper boundary is fixed, and a load of
67,500 N is applied as an evenly distributed pressure force over the
rightmost two ports at the tip of the bracket.

\subsubsection{Optimization results and comparison of aggregation parameters}
\label{sec:bracket_results}
\begin{table}
    \centering
    \begin{tabular}{|c|c||c|c|c|c|c|c|c|}
        \hline
        \(\naggs\) & \(p\) & \(N_{cons}\) & \(N_{jac}\) & Run time (s) & \(m_{frac}^{opt}\) & \(m_{frac}^{pp}\) &
        \(\max \relaxedstress\) (MPa) & \(\max \vonmises\) (MPa) \\
        \hline
        \hline

         1 & 10 & 5114 & 1985 &  9040 & 11.2\% &  8.7\% & 1215 & 938 \\
         5 & 10 &  716 &  394 &  1612 & 14.1\% & 14.1\% & 1148 & 889 \\
        10 & 10 &  903 &  524 &  2345 & 13.0\% & 11.5\% & 1217 & 888 \\
        15 & 10 &  649 &  328 &  1793 & 13.6\% & 13.5\% & 1229 & 884 \\
        20 & 10 &  877 &  457 &  2574 & 13.7\% & 13.3\% & 1215 & 886 \\
        25 & 10 & 1429 &  675 &  4255 & 13.3\% & 12.2\% & 1238 & 886 \\

        \hline
        \hline

         1 & 15 & 4531 & 1495 &  8188 & 15.5\% & 13.9\% &  892 & 907 \\
         5 & 15 & 3348 & 1498 &  6916 & 15.6\% & 13.8\% &  988 & 771 \\
        10 & 15 & 1117 &  536 &  2661 & 15.7\% & 13.7\% &  949 & 799 \\
        15 & 15 & 1113 &  502 &  2832 & 15.5\% & 13.6\% &  929 & 906 \\
        20 & 15 & 1362 &  580 &  3672 & 15.3\% & 12.2\% &  923 & 1074 \\
        25 & 15 & 1046 &  421 &  2935 & 14.8\% & 11.9\% &  957 & 695 \\

        \hline
        \hline

         1 & 25 & 6728 & 2287 & 12580 & 16.4\% & 13.0\% &  747 & 913 \\
         5 & 25 & 1984 & 1097 &  4504 & 17.6\% & 14.4\% &  735 & 704 \\
        10 & 25 & 1162 &  534 &  2835 & 17.9\% & 14.8\% &  725 & 698 \\
        15 & 25 & 1526 &  603 &  3960 & 18.5\% & 14.6\% &  707 & 731 \\
        20 & 25 & 1350 &  574 &  3768 & 18.0\% & 15.2\% &  714 & 686 \\
        25 & 25 & 1377 &  677 &  4167 & 19.5\% & 15.4\% &  725 & 765 \\

        \hline
        \hline

         1 & 50 &    - &    - &     - &      - &      - &    - &   - \\
         5 & 50 & 4228 & 1906 &  8503 & 19.7\% & 17.3\% &  652 & 651 \\
        10 & 50 & 2466 & 1070 &  5701 & 21.2\% & 17.1\% &  610 & 602 \\
        15 & 50 & 2887 & 1233 &  7161 & 21.8\% & 16.5\% &  579 & 617 \\
        20 & 50 & 1984 &  826 &  5442 & 22.1\% & 17.8\% &  557 & 606 \\
        25 & 50 & 2489 & 1068 &  7116 & 22.8\% & 18.0\% &  560 & 592 \\

        \hline
    \end{tabular}
    \caption{Summary of optimization results for varying aggregation parameters.
    For each combination of $\naggs$ and $p$, reported are, from left to right:
    number of constraint evaluations, number of constraint Jacobian evaluations, total
    run time, mass fraction at the local optimum, mass fraction after postprocessing,
    maximum relaxed stress at the optimum, and maximum Von Mises stress after postprocessing.
    Stresses are computed using the full order model in a postprocessing step.}
    \label{tab:agg_comparison}
\end{table}

We solve the optimization \eqref{eq:actual_formulation} for values of the aggregation
multiplier $p =$ 10, 15, 25, and 50, and with 1, 5, 10, 15, 20, and 25 aggregation
regions. The results are summarized in Table~\ref{tab:agg_comparison}, which records
the number of constraint evaluations ($N_{cons}$), constraint Jacobian evaluations
($N_{jac}$), total optimization time, the mass fraction at convergence
($m_{frac}^{opt}$), the mass fraction after postprocessing ($m_{frac}^{pp}$),
the maximum relaxed stress at convergence ($\max \relaxedstress$), and
maximum Von Mises stress in the postprocessed design ($\max \vonmises$).
An optimization not terminating within 10,000
iterations is considered non-convergent, as in the case of the run for $p = 50$
and $\naggs = 1$. For all cases, the maximum stress for the optimization algorithm
is not reduced from the desired maximum: $\heurstressmax = \stressmax$.
All designs are postprocessed using a dropping tolerance of $\densitymin = 0.2$.

\begin{figure}
    \centering
    \includegraphics[width=\textwidth]{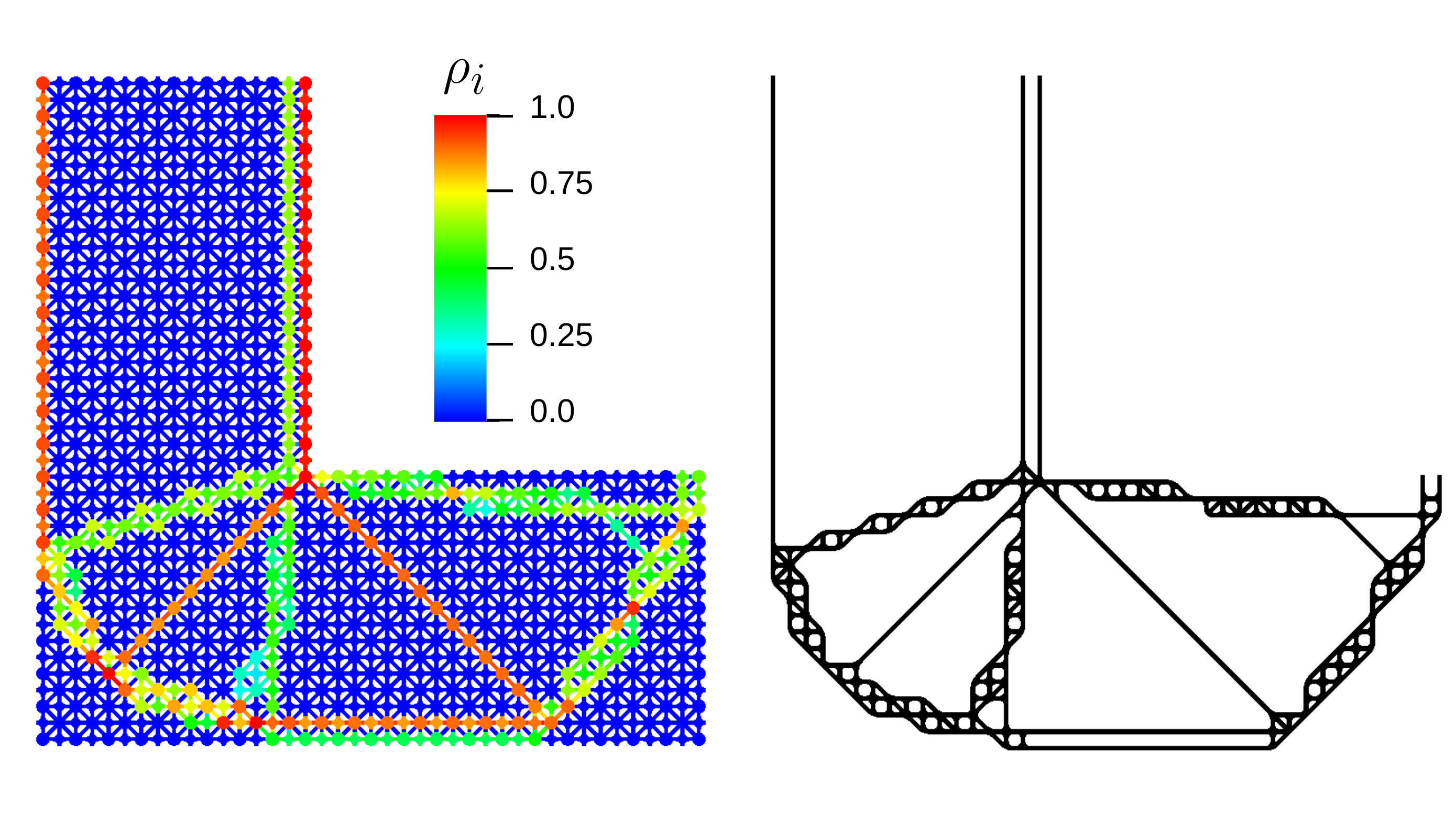}
    \caption{Optimization result using 10 aggregation regions and $p = 15$. Left:
    locally optimal density field. Right: the resulting design after postprocessing.
    Maximum Von Mises stress in the postprocessed design is 799 MPa, compared to
    $\stressmax = 880$ MPa. Before postprocessing, the maximum relaxed stress was
    $\relaxedstress = 949$ MPa.}
    \label{fig:l-bracket-example}
\end{figure}

Fig.~\ref{fig:l-bracket-example} illustrates the result using 10 aggregation
regions and $p = 15$, a representative optimization result. Many components at
the local optimum have intermediate values of density; however, component densities
are clearly divided between material and void components with densities being either
greater than approximately $\rho = 0.25$, or close to 0. This motivates our choice
of dropping tolerance.
Postprocessing creates a fully solid design that satisfies the desired
stress constraint, as well as reducing the mass from the value at the local
optimum, since many components are replaced with lighter substitutes despite
having their density increased to 1 from an intermediate value. The design
uses square lattice cells both with and without diagonal reinforcement in
different areas of the design, as well as long, solid struts and larger open
cells. We highlight here the flexibility of the component-wise approach over
a design method that assumes a functionally graded structure,
with regions composed of periodic unit cells.

Although we observe variability in the optimization's performance,
there are some general trends. First, increasing the aggregation multiplier
$p$ results in the expected behavior: designs are more conservative, and the
optimization problem is more difficult. Increasing the number of aggregation
regions, however, does not have a discernable effect on the value of the maximum
stress attained in the final design, at least for a moderate number of aggregation
regions as used here. It does appear that adding
regions is beneficial for convergence of the algorithm up to 10-15 regions; for values of $p =$ 15, 25,
and 50, the number of constraint evaluations is reduced enough that overall run time
is faster for ten aggregation regions than five, despite requiring twice the number of
adjoint solves for each Jacobian evaluation.

In every case, the postprocessing step results in a mass that is less than or
equal to to the objective value at the local optimum. In a few runs, however,
it does result in violating the stress constraint because the nonlinear
optimization cannot account for the changes made to the shape of components.
Removing components with a near-zero density has a negligible effect on the
structural response, but the component substitution step may have a large effect.
In the examples using one aggregation region and $p = $ 15 and 25, we observe
that component substitution results in increasing the maximum stress to a level
that violates the max constraint, where before substitution it satisfied or
nearly satisfied the constraint. In the example using 20 aggregation regions
and $p =$ 15, a large increase in maximum stress is observed; however, this is
due to the choice of $\densitymin = 0.2$, which is too aggressive for this case
and removes some components that have a smaller density but are structurally
important. For other examples, the postprocessing step either results in a
maximum stress that is less than the maximum relaxed stress at the local
optimum, or increases it by a small amount such that the max constraint is
still respected.

Despite variation across runs, it is clear that with a large enough choice of
aggregation multiplier the proposed methodology consistently converges to
conservative designs with a large mass reduction.

\subsubsection{Performance of the ROM}
\label{sec:rom_comparison}
We compare the accuracy of the ROM used for the primary results to results
using a ROM trained with a 99.99\% criterion for total energy as well as
one using 99\% of total energy. We define the following error measures:
\begin{equation}
    e^{\sigma_r}_{max} = \frac{\max_{\Omega} \tilde{\sigma}_r - \max_{\Omega} \relaxedstress}
    {\max_{\Omega} \relaxedstress}
\end{equation}
\begin{equation}
    e^{\sigma_r} = \frac{\left\|\tilde{\sigma}_r - \relaxedstress\right\|_{L^2} }
    {\left\|\relaxedstress\right\|_{L^2}}
\end{equation}
\begin{equation}
    e^{u} = \frac{\left\| \tilde{u} - u \right\|_{L^2}}{\left\|u\right\|_{L^2}},
\end{equation}
where $\tilde{\sigma}_r$ is the relaxed Von Mises stress as computed from the ROM solution,
$\relaxedstress$ is the relaxed Von Mises stress computed from the full-order solution,
$\tilde{u}$ is the ROM displacement, and $u$ is the full-order displacement. The
$L^2$ norm has its usual definition.

Reported in Table \ref{tab:rom_comparison} are these errors for both the full system
at the locally optimum parameter value, and for the postprocessed design. We report
results for the designs using an aggregation multiplier of 15; the errors reported
for these designs are representative of those for all designs. The exception is error
of 13\% in the postprocessed design for 5 aggregation regions, which is an outlier.
In no other case was the error in maximum stress larger than 10\%. The use of a
99.99\% total energy criterion is sufficient to reduce the error in maximum stress to
less than 5\% in every case, and less than 1\% for all but one; the effect on the
displacement error is smaller. Errors when using
a ROM capturing only 99\% of total energy are large in every case; this criterion
is clearly not strict enough to be useful.

\begin{table}
    \centering
    \begin{tabular}{|c||c|c|c|c|c|c|}
        \hline
        \multicolumn{1}{|c||}{} & \multicolumn{3}{|c|}{Errors at local optimum} & \multicolumn{3}{|c|}{Errors in postprocessed design} \\
        \hline
        $\naggs$ & $e^{\sigma_r}_{max}$ & $e^{\sigma_r}$ & $e^{u}$ & $e^{\sigma_r}_{max}$ & $e^{\sigma_r}$ & $e^u$ \\

        \hline
        \hline

        \multicolumn{7}{|c|}{Errors using 99.9\% of total energy} \\
        \hline
        1 & -0.52\% & 4.12\% & 2.80\% & -0.67\% & 2.96\% & 0.63\% \\
        5 & 3.17\% & 4.46\% & 4.18\% & 13.06\% & 2.93\% & 4.21\% \\
        10 & -0.58\% & 3.67\% & 2.98\% & -0.78\% & 2.84\% & 5.21\% \\
        15 & 6.42\% & 3.91\% & 4.04\% & -0.98\% & 2.80\% & 1.96\% \\
        20 & 0.05\% & 4.17\% & 3.19\% & 0.09\% & 3.22\% & 3.35\% \\
        25 & 2.33\% & 3.78\% & 3.41\% & 8.66\% & 2.54\% & 4.44\% \\

        \hline
        \hline

        \multicolumn{7}{|c|}{Errors using 99.99\% of total energy} \\
        \hline
        1 & 0.54\% & 2.72\% & 1.83\% & 0.02\% & 2.30\% & 0.49\% \\
        5 & 0.08\% & 3.05\% & 2.83\% & -0.90\% & 2.28\% & 4.09\% \\
        10 & 0.02\% & 2.33\% & 2.01\% & 0.09\% & 2.33\% & 5.04\% \\
        15 & 4.03\% & 2.62\% & 2.68\% & -0.08\% & 2.26\% & 1.94\% \\
        20 & -0.79\% & 2.77\% & 2.11\% & -0.18\% & 2.71\% & 3.24\% \\
        25 & -0.98\% & 2.50\% & 2.27\% & 0.24\% & 2.10\% & 4.34\% \\

        \hline
        \hline 

        \multicolumn{7}{|c|}{Errors using 99\% of total energy} \\
        \hline
        1 & 155.67\% & 48.38\% & 25.88\% & 108.61\% & 21.24\% & 8.65\% \\
        5 & 167.64\% & 42.61\% & 29.99\% & 134.32\% & 21.10\% & 47.11\% \\
        10 & 160.05\% & 33.56\% & 21.75\% & 149.29\% & 19.87\% & 55.57\% \\
        15 & 173.93\% & 37.32\% & 30.97\% & 112.66\% & 20.84\% & 23.10\% \\
        20 & 176.95\% & 43.40\% & 23.86\% & 71.62\% & 19.09\% & 38.91\% \\
        25 & 140.72\% & 37.28\% & 25.13\% & 197.87\% & 17.45\% & 46.89\% \\
        \hline
    \end{tabular}
    \caption{Comparison of ROM accuracies for different total energy criteria.
    Reported are relative errors in the maximum stress and relative errors in
    stress and displacement as measured in the $L^2$ norm for both the optimal
    parameter value and the postprocessed design.}
    \label{tab:rom_comparison}
\end{table}

Table~\ref{tab:rom_timings} additionally reports the resulting number of degrees
of freedom and the speedup relative to the full-order static condensation model
for each total energy criterion. We report the speedup for the forward solve,
which includes the cost of factoring the Schur complement matrix and is the largest
component of optimization run time, and the speedup from use of our efficient
aggregation scheme relative to the straightforward method using quadrature. Speedup
is measured for the full system as used during the optimization procedure.
The 99.9\% total energy model achieves a speedup of more than a two
orders of magnitude vs. a full-order static condensation solve, and
the efficient aggregation scheme yields an almost 2.5x speedup vs. straightforward
quadrature. In \cite{our_first_paper}, we found that the full-order static
condensation achieved approximately a 5x speedup compared to the underlying
finite element model; thus, we estimate that the 99.9\% total energy ROM used in
these optimizations will provide approximately a 750x speedup over a finite element
solve using the same mesh.

\begin{table}
    \centering
    \begin{tabular}{|c||c|c|c|}
        \hline
        Criterion & Total DOFs & Forward speedup & Aggregate speedup \\
        \hline
        99\% & 44,214 & 224x & 2.78x \\
        99.9\% & 53,056 & 151x & 2.47x \\
        99.99\% & 70,740 & 59.6x & 1.54x \\
        \hline
    \end{tabular}
    \caption{ROM speedups vs. a full-order static condensation solve for the forward
    solve, and a straightforward quadrature for the computation of stress aggregates.}
    \label{tab:rom_timings}
\end{table}

\subsubsection{Comparison to results using first-order elements}
\label{sec:first_order_comparison}
We also inquire whether the use of second-order finite elements
to construct the ROM is actually beneficial. We expect that the resulting stresses
in element interiors are more accurate, but it is uncertain whether this increased
accuracy has a beneficial effect in the optimization. To determine whether
it does, we ran optimizations using a ROM constructed using first-order elements and
a 99.9\% total energy criterion to compare to the corresponding results using
second-order elements. These optimizations use $p=15$ and 10,
15, 20 and 25 aggregation regions, respectively.
The ROM constructed using first-order elements has exactly the same
dimension as that constructed using second-order elements, and the quadrature rule
used to compute stress aggregates is identical. The aggregation
regions are also identical to those used for the optimizations in
Section \ref{sec:bracket_results}.

\begin{table}
    \centering
    \begin{tabular}{|c||c|c|c|c|c|c|}
        \cline{4-7}
        \multicolumn{3}{c}{} & \multicolumn{2}{|c|}{1st-order elements} &
        \multicolumn{2}{|c|}{2nd-order elements} \\
        \hline
        $n_{agg}$ & Run time (s) & $m_{frac}^{pp}$ & $\max\relaxedstress$ &
        $\max\vonmises$ & $\max\relaxedstress$ & $\max\vonmises$ \\
        \hline
        \hline

        10 & 13936 & 10.6\% & 1063 & 769 & 1215 & 944 \\
        15 &  3287 & 13.5\% & 1020 & 962 & 1195 & 1023 \\
        20 &  2743 & 11.8\% & 1001 & 899 & 1148 & 910 \\
        25 &  4699 & 11.1\% & 1026 & 1027 & 1153 & 1155 \\
        \hline
    \end{tabular}
    \caption{Optimization results using a ROM constructed with first-order elements
    as the underlying model. We show the optimization run time and postprocessed 
    mass fraction from an optimization run using first-order elements,
    and maximum stresses at the optimum and in the postprocessed design
    as computed using both first- and second-order elements. These stresses are
    computed using the corresponding full-order models.}
    \label{tab:first-order-comparison}        
\end{table}

To compare the optimization results we examine the run time, mass fraction of the
postprocessed design, and the max stresses at the optimum and in
the postprocessed design. The latter are computed using both the full-order 
first-order finite element model and the full-order second-order finite element
model to assess the impact of second-order elements on accuracy in the stress.
The results are summarized in Table~\ref{tab:first-order-comparison}; the optimization
run time and mass fraction may be compared to the second block in Table~\ref{tab:agg_comparison}. We find that the first-order model underestimates the maximum
stress significantly. As a result, the optimized mass fractions are smaller than the
corresponding optimization results with second-order elements, but none of the designs
found using first-order elements is conservative when analyzed using second-order elements.
Additionally, in three out of four cases, the optimization with first-order elements as the
model had a longer run time than the second-order optimization. This is due to the fact
that the ROM constructed from first-order elements has the same dimension as that constructed
with second-order elements, while the optimizer spent additional model evaluations on
line search in these cases. These results justify the use of the second-order model
for our primary results.

\subsubsection{Effect of decreasing the maximum stress}
To verify that the heuristic of imposing a reduced stress limit $\heurstressmax$
(Sec.~\ref{sec:conservative_result}) achieves the desired effect, we solve the
optimization problems with $p = 10$ using $\heurstressmax = 800$ MPa. The
corresponding optimizations with $\heurstressmax = 880$ MPa yielded postprocessed
designs that did not satisfy the stress constraint. Results are reported
in Table~\ref{tab:lower_max_results} in the same format as those in Table~\ref{tab:agg_comparison}.

\begin{table}
    \centering
    \begin{tabular}{|c|c||c|c|c|c|c|c|c|}
        \hline
        \(\naggs\) & \(p\) & \(N_{cons}\) & \(N_{jac}\) & Run time (s) & \(m_{frac}^{opt}\) & \(m_{frac}^{pp}\) &
        \(\max \relaxedstress\) (MPa) & \(\max \vonmises\) (MPa) \\
        \hline
        \hline

        \multicolumn{9}{|c|}{$\heurstressmax = 800$ MPa} \\
        \hline
         1 & 10 & 2399 & 1251 &  4572 & 15.2\% & 12.9\% & 1159 & 910 \\
         5 & 10 &  827 &  517 &  1866 & 16.2\% & 14.9\% &  996 & 775 \\
        10 & 10 &  803 &  458 &  2047 & 16.2\% & 14.8\% &  984 & 698 \\
        15 & 10 &  702 &  384 &  1887 & 16.3\% & 15.4\% & 1029 & 704 \\
        20 & 10 &  806 &  484 &  2393 & 15.5\% & 13.6\% &  994 & 885 \\
        25 & 10 &  979 &  564 &  3009 & 16.0\% & 14.1\% &  970 & 906 \\

        \hline
        \hline

        \multicolumn{9}{|c|}{$\heurstressmax = 880$ MPa} \\
        \hline
         1 & 10 & 5114 & 1985 &  9040 & 11.2\% &  8.7\% & 1215 & 938 \\
         5 & 10 &  716 &  394 &  1612 & 14.1\% & 14.1\% & 1148 & 889 \\
        10 & 10 &  903 &  524 &  2345 & 13.0\% & 11.5\% & 1217 & 888 \\
        15 & 10 &  649 &  328 &  1793 & 13.6\% & 13.5\% & 1229 & 884 \\
        20 & 10 &  877 &  457 &  2574 & 13.7\% & 13.3\% & 1215 & 886 \\
        25 & 10 & 1429 &  675 &  4255 & 13.3\% & 12.2\% & 1238 & 886 \\
        \hline
    \end{tabular}
    \caption{Comparison of optimization results when decreasing the value of
    the stress limit for optimization. Consult
    Table~\ref{tab:agg_comparison} for definitions of the notation used here.}
    \label{tab:lower_max_results}
\end{table}

Setting $\heurstressmax < \stressmax$ achieves a reduction in the maximum relaxed
stress at the optimum in every case, as well as accelerating convergence in most
cases. The effect after postprocessing is somewhat unpredictable. In the
case with 25 aggregation regions, the maximum postprocessed stress using a lower
upper bound for optimization is actually higher than that with the higher upper
bound. In all other cases, the postprocessed stress was decreased.
We expect that with enough reduction of $\heurstressmax$, optimization results will
eventually lead to conservative postprocessed designs. Determining how much reduction
is required, however, may require an iterative procedure and a choice that is too
conservative will trade optimality for conservativeness. Instead, it may be preferable
to choose a higher value for the aggregation multiplier $p$, and/or a larger number
of aggregation regions. The choice of these parameters of the optimization is
problem-dependent.

\subsection{Mass minimization of a cantilever beam}

Our second numerical example demonstrates how the component-wise approach 
solves a different structural optimization problem while using the same set of components.
This illustrates a key advantage: the offline phase of constructing a component
library need only be performed once, then the resulting dataset used to solve
multiple design problems.
We minimize the mass of a cantilever beam, fixed at
one end and with a vertical load applied to the opposite tip. The problem setup is
illustrated in Fig.~\ref{fig:beam-setup}; material properties are the same as for
the L-bracket. The cantilever beam has a length of 1 m and height of approximately
25.8 cm (due to the lattice structure); a 30 kN load is applied to the rightmost
two ports on the upper surface of the beam as an evenly distributed pressure force.

\begin{figure}
    \centering
    \includegraphics[width=0.95\textwidth]{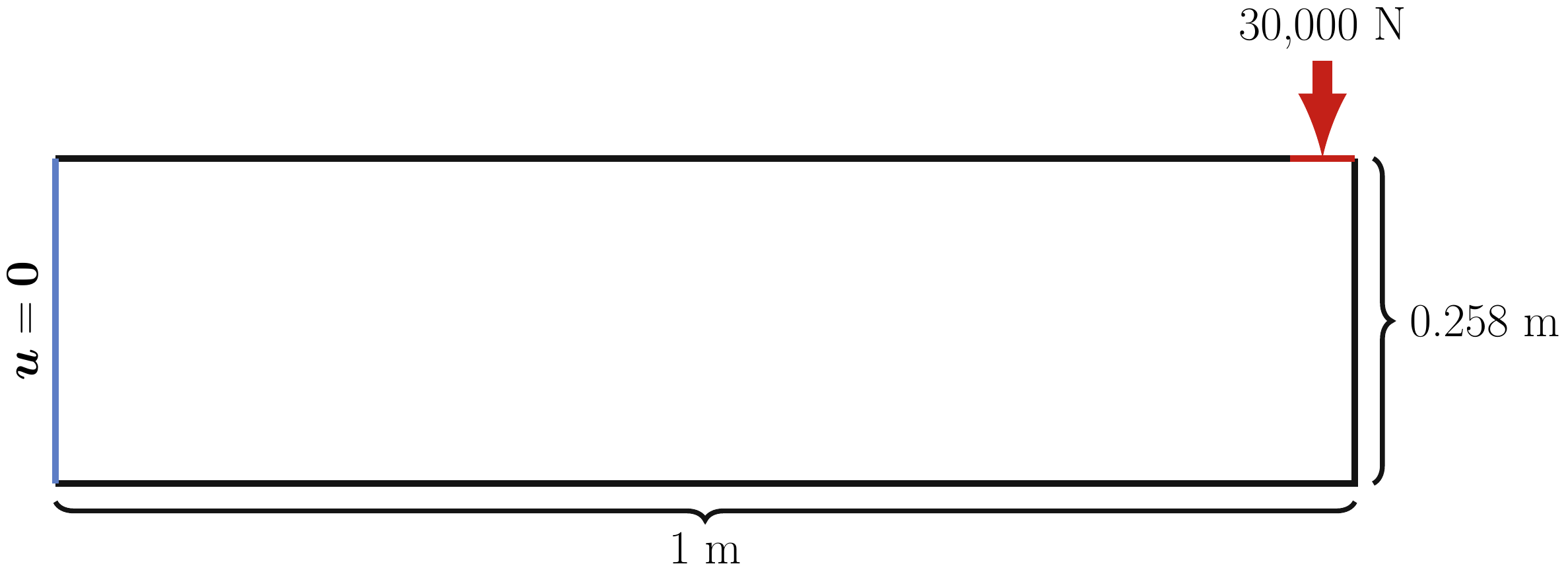}
    \caption{Setup for the cantilever beam optimization.}
    \label{fig:beam-setup}
\end{figure}

The underlying finite element mesh for the beam's ground structure contains
3,750,278 second-order triangular elements for a total of more than 15 million
degrees of freedom (not accounting for those duplicated on internal ports).
The full-order static condensation system has 428,736 degrees of freedom, and
the ROM constructed with 99.9\% of total energy for port spaces contains
81,664.

\subsubsection{Optimization result}
\label{sec:beam_result}
We solve the optimization problem \eqref{eq:actual_formulation} using 10
aggregation regions and an aggregation multiplier of $p = 15$ and the ROM
using a 99.9\% total energy criterion. The convergence tolerance for Ipopt
is set to $10^{-6}$. For this optimization, no reduction of the upper stress
bound $\heurstressmax$ was necessary to obtain a conservative design. The
optimization algorithm converges in 1,944 iterations with 2,576 constraint
evaluations. At the local optimum, the mass fraction relative to the ground
structure is 19.2\%, with a maximum relaxed stress of $\relaxedstress = 849$
MPa. After postprocessing, the design has a mass fraction of 15.1\% and
maximum stress of 784 MPa. The result is postprocessed using dropping tolerance
$\densitymin = 0.25$.

\begin{figure}
    \centering
    \includegraphics[width=0.95\textwidth]{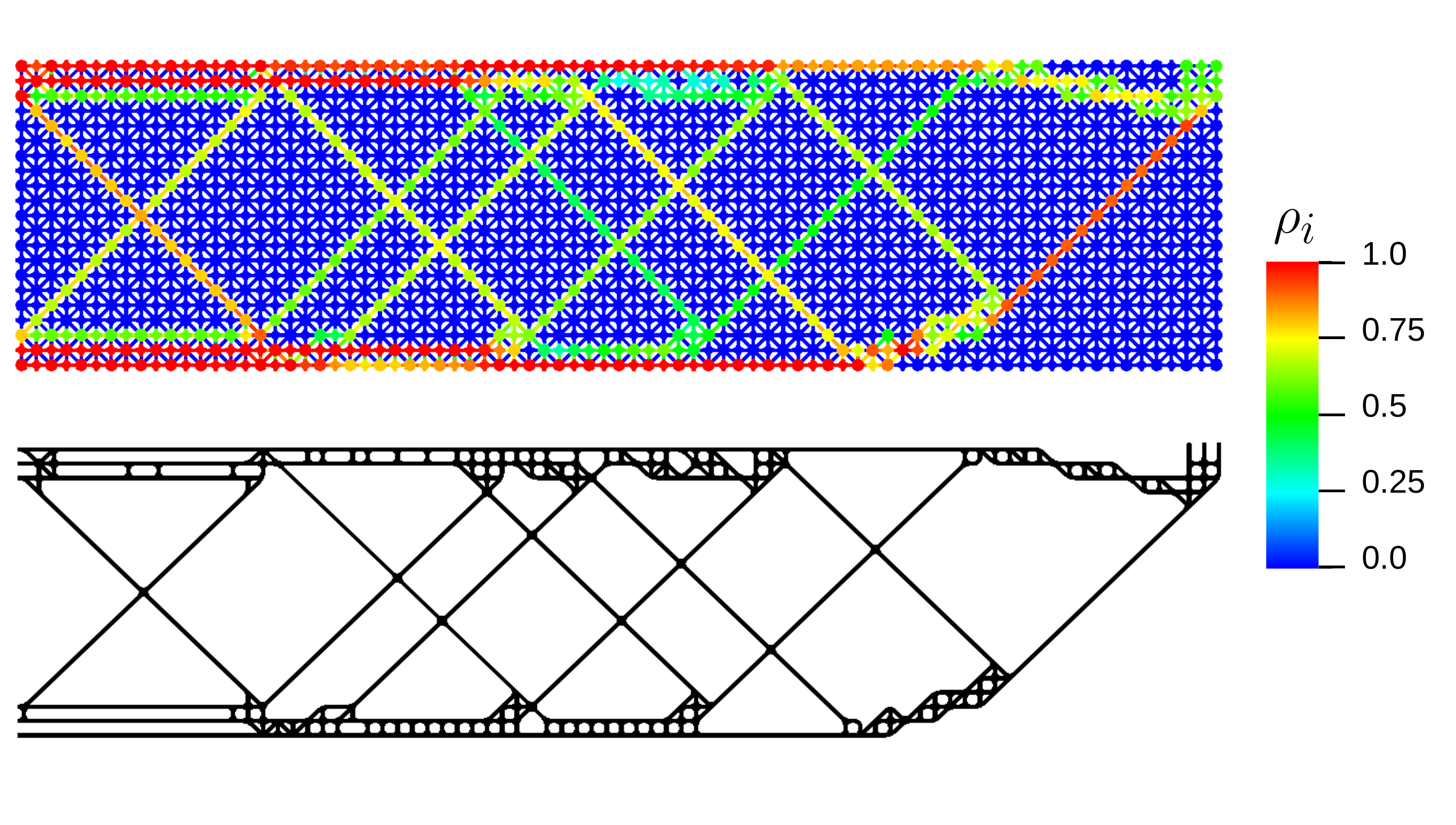}
    \caption{Results of the cantilever beam optimization. Top: optimized
    component-wise density field. Bottom: the postprocessed structure.}
    \label{fig:beam-result-fig}
\end{figure}

Fig.~\ref{fig:beam-result-fig} shows the optimized density field
and the result of postprocessing. As in the previous optimization results,
many components have intermediate densities, but the chosen dropping tolerance
clearly separates material and void regions. The postprocessed design resembles
well-known results for compliance minimization problems and also has
characteristics of functionally-graded structures -- note the use of small
repeating cells on the upper and lower boundaries, and a structure with a
much larger void fraction in the interior.

\subsubsection{Comparison to an optimization using different random aggregation regions}
\label{sec:agg_region_comparison}
To justify the random assignment of aggregation regions, we
compare optimization results for  different choices of aggregation
regions to the result shown previously. The postprocessed designs are compared to
the previous one in Figure \ref{fig:agg-result-comparison}.
The third result is postprocessed using a dropping tolerance of $\densitymin = 0.2$
(vs. $\densitymin = 0.25$ for the first two designs) because of 
component densities at the optimum that fall between these values but
correspond to structurally important components.

\begin{figure}
    \centering
    \includegraphics[width=0.95\textwidth]{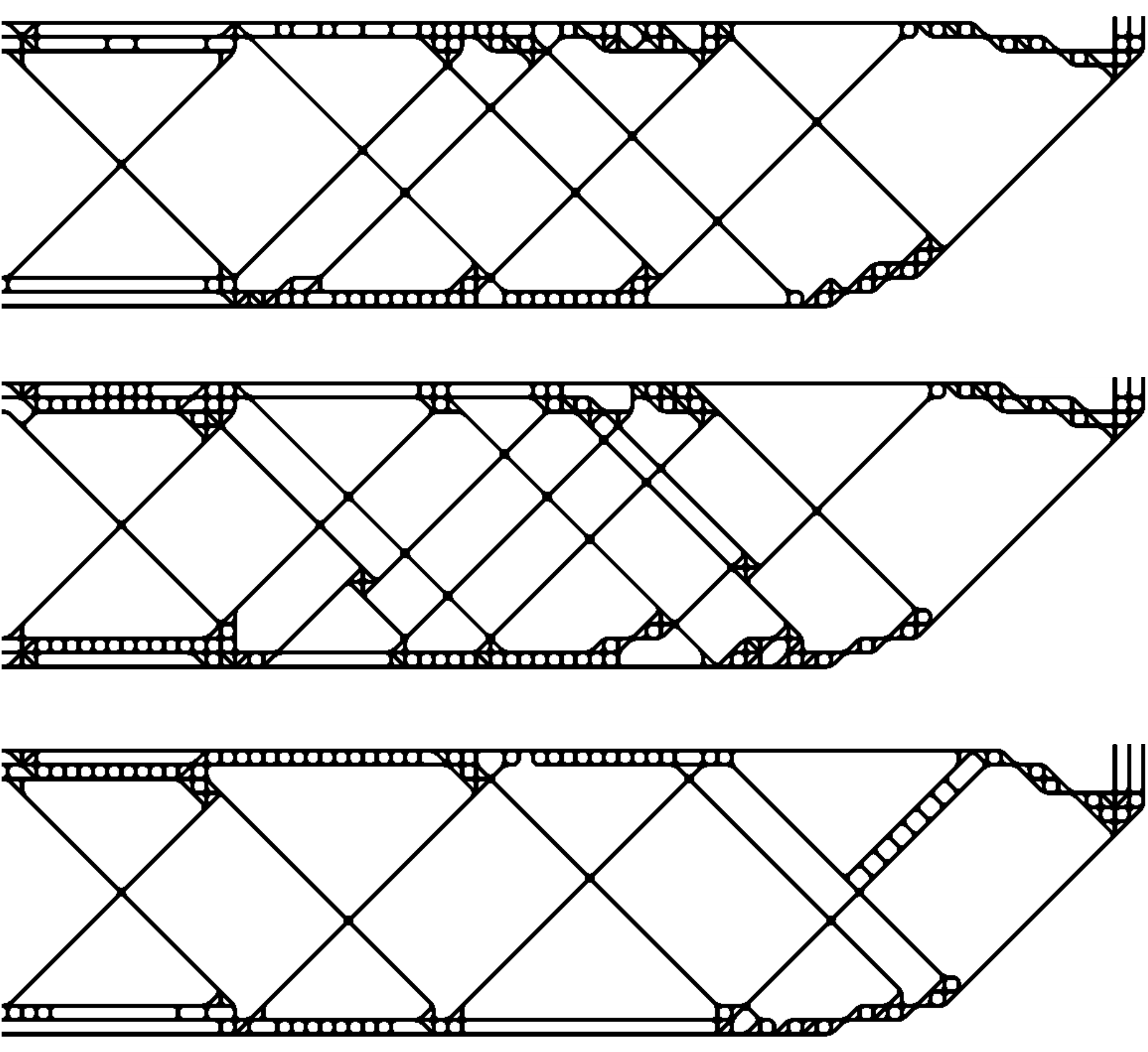}
    \caption{Comparison of the resulting postprocessed designs from
    three different random assignments of aggregation regions. Top: the
    design from Section \ref{sec:beam_result}, with a mass fraction of
    15.1\% and maximum stress of 784 MPa. Middle: the first design with
    redistributed aggregation regions, with mass fraction of 17.4\% and
    a maximum stress of 806 MPa. Bottom: A third design with different
    aggregation regions, with mass fraction of 14.9\% but a maximum stress
    of 1080 MPa.}
    \label{fig:agg-result-comparison}
\end{figure}

Changing aggregation regions does make a significant difference
in optimization results. Whereas the original optimization required
1,944 iterations to converge to tolerance, the optimizations with redistributed
regions required only 793 and 1,428 respectively. We believe this is more
indicative of the difficulty and sensitivity of the optimization problem
than it is of a deficiency in the random strategy for
region assignment.
The local optima with new aggregation regions have mass fractions of 20.3\% and
19\%, vs. 19.2\% in the original; after postprocessing, the mass fractions are
17.4\% and 14.9\% vs. the original 15.1\%. Maximum relaxed stresses at the local
optima are 849 MPa and 887 MPa, vs. the 860 MPa of the original design, and after
postprocessing the new designs have maximum stresses of 806 MPa and 1080 MPa.
Clearly the assignment of aggregation regions influences the quality of optimization
results. We note, however, that the third design's violation of the maximum stress
constraint is due to postprocessing, not the assignment of aggregation regions;
postprocessing creates a stress concentration that does not exist at the local optimum.
Otherwise, random region assignments do result in fairly consistent designs.

\section{Conclusions}
We have presented a novel application of component-wise ROMs to stress constrained
TO. The methodology succeeds in finding solutions with a large mass reduction relative
to the ground structure while respecting stress constraints.
While the presence of many local optima as well as the approximation of the
maximum stress mean that it is impossible to
guarantee a conservative solution, adjusting the aggregation multiplier and
upper bound on stress proved effective for stress control.
The reduced order model provides an error in the maximum stress
of less than 10\% in most cases and a relative error of about 3\% in the $L^2$ norm,
while providing a 150x speedup relative to a full order static condensation model.
If greater accuracy is required, a ROM using 99.99\% of total energy still provides
a 60x speedup while reducing the error in the maximum stress by approximately an
order of magnitude. Future work may investigate adaptive refinement of the ROMs to
use the more accurate model only for components where it is needed, as
permitted by the component-wise discretization. Due to the acceleration provided
by the component-wise ROM, we are able to solve the TO problem efficiently on
a ground structure containing millions of second-order finite elements. Component
ROMs provide a high accuracy surrogate model that may be reused to study a
multitude of both homogeneous and heterogeneous lattice-like structures.

\section*{Declaration of competing interest}
The authors declare that they have no known competing financial interests or personal relationships that could
have appeared to influence the work reported in this paper.

\section*{Acknowledgements}
Funding: This work was supported in part by the AEOLUS center under US Department of
Energy Applied Mathematics MMICC award DE-SC0019303. The second author was funded by
LDRD (21-FS-042) at Lawrence Livermore National Laboratory. Lawrence Livermore
National Laboratory is operated by Lawrence Livermore National Security, LLC,
for the U.S. Department of Energy, National Nuclear Security Administration
under Contract DE-AC52-07NA27344 and LLNL-JRNL-834458.

\printbibliography

\end{document}